\documentclass[11pt,reqno]{amsart}
\usepackage[T1]{fontenc}
\usepackage[utf8]{inputenc}
\usepackage{lmodern}
\usepackage[a4paper,margin=1.08in]{geometry}
\usepackage{amsmath,amssymb,amsfonts,amsthm,mathtools,mathrsfs}
\usepackage{microtype,enumitem,url,needspace,xcolor}
\usepackage[hidelinks]{hyperref}
\allowdisplaybreaks
\numberwithin{equation}{section}
\newtheorem{theorem}{Theorem}[section]
\newtheorem{lemma}[theorem]{Lemma}
\newcommand{\M}{\mathcal M}
\newcommand{\A}{\mathcal A}
\newcommand{\Q}{\mathcal Q}
\newcommand{\R}{\mathbb R}
\newcommand{\Z}{\mathbb Z}
\newcommand{\C}{\mathbb C}
\newcommand{\T}{\mathbb T}
\newcommand{\E}{\mathbb E}
\newcommand{\one}{\mathbf 1}
\newcommand{\dd}{\,\mathrm d}
\newcommand{\bd}{b_{\mathrm d}}
\newcommand{\boff}{b_{\mathrm{off}}}
\newcommand{\lo}{\mathrm{lo}}
\newcommand{\hi}{\mathrm{hi}}
\theoremstyle{remark}
\newtheorem{remark}[theorem]{Remark}
\theoremstyle{remark}
\newcommand{\core}{\mathscr D_{\M}}
\newcommand{\norm}[2][]{\lVert #2\rVert_{#1}}
\DeclareMathOperator{\supp}{supp}
\DeclareMathOperator{\Id}{Id}
\DeclareMathOperator{\pv}{p.v.}
\makeatletter
\@namedef{subjclassname@2020}{\textup{2020} Mathematics Subject Classification}
\makeatother
\title[Noncommutative maximal rough singular integrals]
{Weak type $(1,1)$ estimates for noncommutative\\
	maximal truncated rough singular integrals}
\date{}
\author[Cao]{Junming Cao}
\address{Junming Cao, School of Mathematical Sciences, Huaqiao University, Quanzhou 362021, China} \email{caojunming@hqu.edu.cn}

\author[Li]{Wei Li}
\address{Wei Li, School of Mathematics and Statistics, Central South University, Changsha 410075, China}
\email{liwei.csu@csu.edu.cn}

\makeatletter
\@namedef{subjclassname@2020}{%
	\textup{2020} Mathematics Subject Classification}
\makeatother

\hypersetup{pdftitle={Weak type (1,1) estimates for noncommutative maximal truncated rough singular integrals}}
\begin{document}
	\begin{abstract}
		We prove a weak type $(1,1)$ estimate for maximal truncations of noncommutative rough  singular integrals whose angular kernels satisfy the classical $L\log L$ integrability and cancellation conditions. The main ingredient of the proof is a new linearization method for handling the maximal truncations acting on the bad functions arising from the noncommutative Calder\'on--Zygmund decomposition.
	\end{abstract}
	\maketitle

\section{Introduction}

Let $d\ge2$ and let $\Omega\in L_1(S^{d-1})$ be
complex-valued with mean zero. The associated rough
homogeneous singular integral is defined, initially
for smooth compactly supported functions, by
 $$
	T_\Omega f(x)=\pv\int_{\R^d}
	\frac{\Omega(y/|y|)}{|y|^d}f(x-y)\dd y.
$$
The weak type $(1,1)$ estimate was proved in dimension
two by Christ and Rubio de Francia \cite{ChristRubio},
and subsequently in all dimensions by Seeger
\cite{Seeger}, under the cancellation and angular
$L\log L$ conditions. Seeger's proof is based on a
microlocal decomposition of the annular kernels.

The corresponding maximal truncated operator is
defined by
\[
T_{\Omega,*}f(x)
=\sup_{\varepsilon>0}
\left|
\int_{|y|>\varepsilon}
\frac{\Omega(y/|y|)}{|y|^d}f(x-y)\dd y
\right|.
\]
Calder\'on and Zygmund \cite{CZ1956} introduced this
operator and established its strong $L_p$ boundedness
for $1<p<\infty$ by the method of rotations.
Its weak type $(1,1)$ estimate under the same
cancellation and angular $L\log L$ conditions
remained open for decades and was recently
established by Lai \cite{LaiMax}.

In the noncommutative setting, Parcet \cite{Parcet}
introduced the noncommutative Calder\'on--Zygmund
decomposition based on Cuculescu projections.
Its refinement by Cadilhac, Conde-Alonso and
Parcet \cite{CCP} separates the bounded good
function from the diagonal and off-diagonal
bad functions. Hong, Lai and Xu \cite{HLX}
proved maximal strong type $(p,p)$ inequalities
for noncommutative  rough singular integrals and maximal weak
type $(1,1)$ bounds under stronger kernel
regularity. The corresponding noncommutative maximal weak type $(1,1)$ estimate
under only the cancellation and angular $L\log L$
conditions remained open. 
 In this paper, we 
  resolve  this problem.

Let $(\M,\tau)$ be a semifinite von Neumann algebra
with a normal semifinite faithful trace, and put
$
\A=L_\infty(\R^d)\,\overline\otimes\,\M,
$ $
\varphi=\int_{\R^d}\dd x\otimes\tau.
$ 
Assume that $\Omega$ satisfies the cancellation
condition
\begin{equation}\label{eq:meanzero}
	\int_{S^{d-1}}\Omega(\theta)\dd\sigma(\theta)=0,
\end{equation}
where $\sigma$ denotes surface measure. Write
$
K_\Omega(y)=\frac{\Omega(y/|y|)}{|y|^d}
$
and define the truncated operators
\[
T_{\Omega,\varepsilon}f(x)
=\int_{|y|>\varepsilon}
K_\Omega(y)f(x-y)\dd y,
\qquad \varepsilon>0,
\]
initially for
$f\in\core=C_c^\infty(\R^d)
\otimes(\M\cap L_1(\M))$.For $\Omega\ne0$, set
\[
C_{\Omega,1}
=\int_{S^{d-1}}|\Omega(\theta)|
\left(1+\log^+
\frac{|\Omega(\theta)|}{\|\Omega\|_1}\right)
\dd\sigma(\theta),
\]
where $\|\Omega\|_1$ denotes the
$L_1(S^{d-1})$ norm; put $C_{0,1}=0$.
Its finiteness is
equivalent to
$\Omega\in L(\log^+L)(S^{d-1})$.

Our main result is stated as below.
\begin{theorem}\label{thm:main}
	Suppose that $\Omega\in L(\log^+L)(S^{d-1})$
	satisfies \eqref{eq:meanzero}. For every
	$f\in L_1(\A)$ and $\lambda>0$, there is a
	projection $e\in\A$ such that
$$
		\sup_{\varepsilon>0}
		\|eT_{\Omega,\varepsilon}f e\|_\infty
		 \le\lambda, \quad
		\lambda\varphi(1-e)
		 \le C_d C_{\Omega,1}\|f\|_1. 
$$
	The constant $C_d$ depends only on the dimension.
\end{theorem}
\begin{remark}
	For $\M=\C$, Theorem~\ref{thm:main} is the  weak type $(1,1)$
	maximal estimate of Lai \cite{LaiMax}.  We also extend the result to the
	nonhomogeneous convolution kernels; see Section~\ref{sec:nonhomogeneous}.
\end{remark}

	The main difficulties in proving Theorem~\ref{thm:main}  arise from two aspects: the roughness of the kernel and the maximal nature of the operator. The former difficulty is already present in the nonmaximal case, where the authors of \cite{CaoLiWu} adapted Seeger's microlocal method to the bad functions arising from the noncommutative Calder\'on--Zygmund decomposition and established the weak type $(1,1)$ estimate. The main challenge in the present paper is therefore to handle the additional difficulty caused by maximal truncations.
		In the scalar case, Lai \cite{LaiMax} used an iterative partition of dyadic cubes to reorganize the localized terms into ordered partial sums, allowing the Rademacher--Menshov theorem to reduce the required maximal $L_2$ estimates to $L_2$ estimates for signed sums.
		This approach does not extend directly to the noncommutative setting, where the Calder\'on--Zygmund decomposition produces additional off-diagonal terms that lack the global $L_1$ summability available in the scalar case. We therefore develop a different linearization method to handle the maximal truncations acting on the bad functions. After decomposing the kernel into dyadic annuli, the main task is to control the corresponding convolution tails uniformly in the truncation parameter. Unlike the iterative partition in \cite{LaiMax}, our construction uses nested enlarged dyadic grids to obtain suitable support projections for the bad terms. Ordered products of these projections then reduce the maximal estimates to the required $L_2$ estimates for suitably transformed full sums. This construction also allows us to adapt the microlocal method of \cite{CaoLiWu} to these sums, for both the diagonal and off-diagonal bad terms.
	
The paper is organized as follows. Section~\ref{sec:pre} introduces the notation, the noncommutative Calder\'on--Zygmund decomposition, and the kernel decompositions used throughout the paper. Section~\ref{sec:main-proof} proves Theorem~\ref{thm:main} assuming several auxiliary lemmas, whose proofs are given in the subsequent sections. Section~\ref{sec:products} proves the lemma on ordered products of projections. Section~\ref{sec:geometry} establishes the dyadic grid construction and the corresponding lemmas. Section~\ref{sec:microlocal} proves the maximal estimates for the diagonal and off-diagonal bad terms. Section~\ref{sec:off-good} proves the maximal $L_2$ estimate for the off-diagonal good terms. Appendix~\ref{sec:nonhomogeneous} extends Theorem~\ref{thm:main} to nonhomogeneous convolution kernels.
	Throughout, $C_d$ denotes a positive constant depending only on $d$;
	its value may change from line to line.
	
	\section{Preliminaries}\label{sec:pre}
	
\subsection{Noncommutative measure spaces}

Let $(\M,\tau)$ be the semifinite von Neumann algebra introduced
above. Denote by $\mathcal P(\M)$ its projections and by
$\mathcal S(\M)_+$ the positive operators in $\M$ whose support
projections have finite trace. Let $\mathcal S(\M)$ be the linear
span of $\mathcal S(\M)_+$. For $1\le p<\infty$, the space
$L_p(\M)$ is the completion of $\mathcal S(\M)$ under
\[
\|x\|_p=\tau(|x|^p)^{1/p},\qquad |x|=(x^*x)^{1/2},
\]
while $L_\infty(\M)=\M$ with the operator norm.

A closed densely defined operator affiliated with $\M$ is
$\tau$-measurable if $\tau(\one_{(t,\infty)}(|x|))<\infty$
for some $t>0$. Write $L_0(\M)$ for the algebra of such operators.
Then
\[
L_p(\M)=\{x\in L_0(\M):\tau(|x|^p)<\infty\},
\qquad 1\le p<\infty.
\]
For $x\in L_0(\M)$, define
\[
d_x(t)=\tau(\one_{(t,\infty)}(|x|)),\qquad
\|x\|_{p,\infty}=\sup_{t>0}t\,d_x(t)^{1/p}.
\]
The weak space $L_{p,\infty}(\M)$ consists of the operators for
which this quasi-norm is finite.

For the algebra $\A=L_\infty(\R^d)\,\overline\otimes\,\M$
with trace $\varphi$, we use the identification
\[
L_p(\A)=L_p(\R^d;L_p(\M)),\qquad
\|f\|_p^p=\int_{\R^d}\tau(|f(x)|^p)\dd x.
\]
All the preceding definitions apply to $(\A,\varphi)$.
For a positive operator $a$, $\supp(a)$ denotes its support
projection. 

\subsubsection{Vector-valued noncommutative $L_p$-spaces}
Let $\Lambda$ be an arbitrary index set. For $1\le p\le\infty$,
a family $(x_\lambda)_{\lambda\in\Lambda}$ belongs to
$L_p(\A;\ell_\infty(\Lambda))$ if it has a factorization
\[
x_\lambda=ay_\lambda b,\qquad
a,b\in L_{2p}(\A),\quad
\sup_{\lambda\in\Lambda}\|y_\lambda\|_\infty<\infty,
\]
where $a,b$ do not depend on $\lambda$. Its norm is
\[
\|(x_\lambda)\|_{L_p(\A;\ell_\infty(\Lambda))}
=\inf_{x_\lambda=ay_\lambda b}
\|a\|_{2p}\sup_{\lambda\in\Lambda}\|y_\lambda\|_\infty\|b\|_{2p}.
\]
For $p=\infty$, this is $\sup_\lambda\|x_\lambda\|_\infty$.
The description  of maximal $L_p$-spaces is  standard; see
\cite{JX,Mei,HLX}.

For $1\le p<\infty$, the weak maximal quasi-norm is
\[
 	 \|(x_\lambda)\|_{\Lambda_{p,\infty}(\A;\ell_\infty(\Lambda))}
  =\sup_{t>0}t\inf\big\{\varphi(1-e)^{1/p}:
	e\in\mathcal P(\A),\
	\sup_{\lambda\in\Lambda}\|ex_\lambda e\|_\infty\le t\big\}.
 \]
A family $(T_\lambda)_\lambda$ is said to be of weak type weak type $(1,1)$ if, for every
$f\in L_1(\A)$ and $t>0$, one can find a projection $e\in\A$
such that
\[
\sup_\lambda\|eT_\lambda(f)e\|_\infty\le t,
\qquad t\varphi(1-e)\le C\|f\|_1.
\]

\subsubsection{Maximal inequalities}

We use the noncommutative Hardy--Littlewood maximal inequality
in the above form. For $r>0$, let
\[
A_rh=|B(0,r)|^{-1}\int_{|y|<r}h(\cdot-y)\dd y.
\]
\begin{lemma}[\cite{Mei}]\label{lem:HL-input}
	For $1<p<\infty$, for any $h\in L_p(\A),$ we have 
	\[
		\|(A_rh)_{r>0}\|_{L_p(\A;\ell_\infty)}
	\le C_{d,p}\|h\|_p.
		 	\]
\end{lemma} By  duality, this gives
\begin{equation}\label{eq:dual-HL}
	\left\|\sum_i A_{r_i}a_i\right\|_p
	\le C_{d,p}\left\|\sum_i a_i\right\|_p,
	a_i\ge0.
\end{equation}

The following estimate is basic.
\begin{lemma}\label{lem:L2-cut}
	If $(x_i)\in L_2(\A;\ell_\infty)$ and $t>0$, there is a
	projection $e\in\A$ such that
	\[
	\sup_i\|ex_ie\|_\infty\le t,\qquad
	\varphi(1-e)\le8t^{-2}
	\|(x_i)\|_{L_2(\A;\ell_\infty)}^2.
	\]
\end{lemma}

We also use the following estimate \cite[Section~2]{JX}.
\begin{lemma}\label{lem:absolute-L1}
	If $\sum_i\|v_i\|_1<\infty$ and $t>0$, there is a projection
	$e\in\A$ such that
	\[
	\varphi(1-e)\le\frac2t\sum_i\|v_i\|_1,\qquad
	\left\|e\left(\sum_i c_iv_i\right)e\right\|_\infty\le t
	\quad (|c_i|\le1).
	\]
\end{lemma}
	\subsection{The dyadic martingales}
	
	For $k\in\Z$, let $\Q_k$ be the family of standard dyadic cubes in
	$\R^d$ of side length $2^{-k}$, and let $\Sigma_k$ be the
	$\sigma$-algebra generated by $\Q_k$. Set
	\[
	\A_k=L_\infty(\R^d,\Sigma_k)\,\overline\otimes\,\M.
	\]
	Then $(\A_k)_{k\in\Z}$ is an increasing filtration of $\A$. Let
	\[
	\E_k:\A\longrightarrow\A_k
	\]
	be the trace-preserving conditional expectation, extended
	contractively to $L_p(\A)$ for $1\leq p<\infty$. For a locally
	integrable operator-valued function $f$,
	\[
	\E_kf=\sum_{Q\in\Q_k}f_Q\one_Q,
	\qquad
	f_Q=\frac1{|Q|}\int_Qf(y)\dd y.
	\]
	We write
	\[
	f_k=\E_kf,
	\qquad
	df_k=f_k-f_{k-1}.
	\]
	The conditional expectations satisfy the following standard properties.
	\begin{enumerate}[label=\textup{(\roman*)},leftmargin=2.2em]
		\item For $k,l\in\Z$,
		\[
		\E_k\E_l=\E_{\min\{k,l\}}.
		\]
		
		\item If $a,b\in\A_k$ and $h\in L_p(\A)$, $1\leq p<\infty$, then
		\[
		\E_k(ahb)=a\,\E_k(h)\,b.
		\]
	\end{enumerate}
 \subsection{Noncommutative Calder\'on--Zygmund decomposition}
 
 In this subsection, we first work under the assumptions
 \begin{equation}\label{eq:finitemodel}
 	f\in L_1(\A)_+\cap\A,
 	\qquad
 	\supp_{\R^d}f\ \text{is compact},
 	\qquad
 	f=\E_Mf\quad\text{for some }M\in\Z.
 \end{equation}
 We begin with Cuculescu's construction; see
 \cite[Lemma~3.1]{Parcet} and \cite[Section~1]{CCP}.
 
 \begin{lemma}\label{lem:cuculescu}
 	Let $f$ satisfy \eqref{eq:finitemodel}	and  $\alpha>0$. There exists a decreasing sequence of projections
 	$(q_k)_{k\in\Z}$ such that $q_k=1$ for all sufficiently negative
 	$k$, and, for the remaining $k$,
 	\[
 	q_k
 	=
 	\one_{[0,\alpha]}(q_{k-1}f_kq_{k-1})
 	\quad\text{in }q_{k-1}\A q_{k-1}.
 	\]
 	Moreover, for every $k\in\Z$,
 	\begin{enumerate}[label=\textup{(\roman*)},leftmargin=2.2em]
 		\item $q_k\in\A_k$ and $
 		q_kf_kq_k\leq\alpha q_k.
 		$
 		\item $q_k$ commutes with $q_{k-1}f_kq_{k-1}$;
 		\item The following  assertion holds
 		$$\alpha\varphi(1-q)\leq\norm[1]{f},$$
 		where
 		$q=\bigwedge_{k\in\Z}q_k.$
 	\end{enumerate}
 \end{lemma}
 Set
 \[
 p_k=q_{k-1}-q_k.
 \]
 Since $q_k\in\A_k$, we may write
 \[
 q_k=\sum_{Q\in\Q_k}q_Q\one_Q,
 \]
 where each $q_Q$ is a projection in $\M$. If $\widehat Q$ denotes
 the dyadic parent of $Q$, define
 \[
 \pi_Q=q_{\widehat Q}-q_Q,
 \qquad Q\in\Q_k.
 \]
 Then
 \[
 p_k=q_{k-1}-q_k
 =\sum_{Q\in\Q_k}\pi_Q\one_Q.
 \]
 The following standard properties of the projections $(p_k)_{k\in\Z}$
 follow from Cuculescu's construction and the dyadic regularity of the filtration, which can be found in
 \cite{Parcet}.
 
 \begin{lemma}\label{lem:pk-properties}
 	Let $f$ satisfy \eqref{eq:finitemodel}.	For the projections defined above, the following assertions hold.
 	\begin{enumerate}[label=\textup{(\roman*)},leftmargin=2.2em]
 		\item
 		For every $k\in\Z$,
 		$
 		p_k\in\A_k.
 		$
 		Moreover, the projections $(p_k)_{k\in\Z}$ are pairwise
 		orthogonal and
 		\[
 		\sum_{k\in\Z}p_k=1-q.
 		\]
 		
 		\item
 		For every $k\in\Z$,
 		\[
 		p_kf_kq_k=q_kf_kp_k=0.
 		\]
 		
 		\item
 		We have
 		\[
 			\alpha\varphi(p_k)\leq\varphi(fp_k),
 			\qquad
 			\sum_{k\in\Z}\varphi(fp_k)\leq\norm[1]{f}.
 		\]
 		
 		\item
 		The dyadic regularity gives
 		\[
 			p_kf_kp_k\leq2^d\alpha p_k,
 			\qquad
 			\pi_Qf_Q\pi_Q\leq2^d\alpha\pi_Q
 			\quad(Q\in\Q_k).
 		\]
 	\end{enumerate}
 \end{lemma}
 
 The noncommutative Calder\'on--Zygmund decomposition of
 \cite[(1.1) and Lemma~1.1]{CCP}, also stated in
 \cite[Theorem~2.5]{HLX}, takes the following form.

 \begin{theorem}\label{thm:regularcz}
 	Let $f$ satisfy \eqref{eq:finitemodel} and let $\alpha>0$.
 	Let $(q_k)_k$, $(p_k)_k$, and $q$ be the projections introduced above.
 	Then
 	\[
 		f=g+\bd+\boff,
 	\]
 	where the following assertions hold.
 	\begin{enumerate}[label=\textup{(\roman*)},leftmargin=2.2em]
 		\item
 		$
 		g=qfq+\sum_kp_kf_kp_k,
 		$
 		with
 		\begin{equation}\label{eq:czbounds}
 			\norm[1]{g}\leq\norm[1]{f},
 			\qquad
 			\norm[\infty]{g}\leq2^d\alpha.
 		\end{equation}
 		
 		\item
 		$
 		\bd=\sum_kb_{\mathrm d,k},
 		$  where $$
 		b_{\mathrm d,k}=p_k(f-f_k)p_k,
 		$$
 		satisfying
 		\[
 		\E_kb_{\mathrm d,k}=0,
 		\qquad
 		\sum_k\norm[1]{b_{\mathrm d,k}}
 		\leq2\norm[1]{f}.
 		\]
 		
 		\item
 		$
 		\boff=\sum_kb_{\mathrm{off},k},
 		$
 		where
 		\[
 		b_{\mathrm{off},k}
 		=p_k(f-f_k)q_k+q_k(f-f_k)p_k
 		=p_kfq_k+q_kfp_k.
 		\]
 	\end{enumerate}
 \end{theorem}
 
 Fix $n_0=100$.

 The cube decomposition of the diagonal bad part and the associated
 localization projection will be used throughout the proof.
 
 \begin{lemma}\label{lem:cz-localization}
 	For $k\in\Z$ and $Q\in\Q_k$, set
 	\[
 	b_Q=\one_Q\pi_Q(f-f_Q)\pi_Q,
 	\]
 	and define
 	\[
 	\zeta_0
 	=\left(
 	\bigvee_k\bigvee_{Q\in\Q_k}
 	\pi_Q\one_{2^{n_0+10}Q}
 	\right)^\perp.
 	\]
 	Then the following assertions hold.
 	\begin{enumerate}[label=\textup{(\roman*)},leftmargin=2.2em]
 		\item For every $k\in\Z$,
 		\[
 		b_{\mathrm d,k}=\sum_{Q\in\Q_k}b_Q,
 		\qquad
 		\int_Qb_Q=0.
 		\]
 		Moreover,
 		\[
 		\sum_{Q\in\Q_k}\norm[1]{b_Q}
 		=\norm[1]{b_{\mathrm d,k}}.
 		\]
 		
 		\item The localization projection satisfies
 		\[
 			\varphi(1-\zeta_0)
 			\leq C_d\alpha^{-1}\norm[1]{f}.
 		\]
 		
 		\item For every $Q\in\Q_k$,
 		\[
 		\zeta_0(x)\pi_Q=\pi_Q\zeta_0(x)=0,
 		\qquad x\in2^{n_0+10}Q.
 		\]
 	\end{enumerate}
 \end{lemma}

 We next record the localization consequence used below; compare
 \cite[(3-3)]{Lai}.
 
 \begin{lemma}\label{lem:local-support}
 	Let $k,j\in\Z$ and $h_k\in L_1(\A)$ satisfy
 	\[
 	(1-p_k)h_k(1-p_k)=0.
 	\]
 	Let $J_j\in L_1(\R^d)$ be a scalar kernel such that
 	\[
 	\supp J_j\subset\{x:|x|\leq2^{j+1}\}.
 	\]
 	Then
 	\begin{equation}\label{eq:localsupport}
 		\zeta_0(J_j*h_k)\zeta_0=0
 		\qquad(j+k<n_0).
 	\end{equation}
 \end{lemma}
 
\subsection{The functions
	\texorpdfstring{$b_{k,\ell}$ and $g_{k,\ell}$}{b and g}}\label{subsec-b-g}

For the analysis of the off-diagonal bad term $\boff$ in
Theorem~\ref{thm:regularcz}, we use two families introduced in
Parcet's original decomposition \cite{Parcet}.	
 For $k\in\Z$ and $\ell\ge1$, define
 \[
 \begin{aligned}
 b_{k,\ell}
 &=p_k(f-f_{k+\ell})p_{k+\ell}
 +p_{k+\ell}(f-f_{k+\ell})p_k,\\
 g_{k,\ell}
 &=p_kdf_{k+\ell}q_{k+\ell-1}
 +q_{k+\ell-1}df_{k+\ell}p_k.
 \end{aligned}
 \]
 By \cite[Lemma~2.13]{CaoLiWu}, for every $\ell\ge1$,
 \begin{equation}\label{eq:lag-block-properties}
 \begin{gathered}
 \sum_k\|b_{k,\ell}\|_1\lesssim_d\|f\|_1,
 \qquad \E_{k+\ell}b_{k,\ell}=0,\\
 \sum_k\|g_{k,\ell}\|_2^2\lesssim_d\alpha\|f\|_1,
 \qquad g_{k,\ell}\in\A_{k+\ell},\\
 \E_{k+\ell-1}g_{k,\ell}=0.
 \end{gathered}
 \end{equation}
 We also use
 $\sum_k\|g_{k,\ell}\|_1\lesssim_d\|f\|_1$.
 Under \eqref{eq:finitemodel},
 \[
 \boff=\sum_k\sum_{\ell\ge1}
 (b_{k,\ell}+g_{k,\ell}).
 \]

For integers $s\ge n_0$, retain the combined functions of
	\cite[Section~6]{CaoLiWu}:
	\begin{equation}\label{eq:combined-block}
		u_j^s=\sum_{s\le n<2s}b_{\mathrm d,n-j}
		+\sum_{\substack{n\ge n_0,\ \ell\ge1\\s\le n+\ell<2s}}
		(b_{n-j,\ell}+g_{n-j,\ell}).
	\end{equation}
	\begin{lemma}[\cite{CaoLiWu}, Lemma~3.2 and \textup{(6.4)}]
		\label{lem:combined-input}
		For every integer $s\ge n_0$,
		\[
		\sum_j\|u_j^s\|_1\lesssim_d s\|f\|_1.
		\]
	\end{lemma}
	
	\subsection{Kernels and the microlocal estimates}
	
	Use the radial partition and the microlocal decomposition from
	\cite[Section~2.3]{CaoLiWu}. Thus
	$\beta\in C_c^\infty((1/2,2))$ is nonnegative and
	$\sum_j\beta(2^{-j}t)^2=1$ for $t>0$, and
	\begin{equation}\label{eq:dyadickernels}
		K_j^\Psi(x)=\beta(2^{-j}|x|)^2
		\frac{\Psi(x/|x|)}{|x|^d},
		\qquad T_j^\Psi h=K_j^\Psi*h.
	\end{equation}
	We use the Fourier transform
	$\widehat h(\xi)=\int e^{-i\langle x,\xi\rangle}h(x)\dd x$.
	The kernel bounds \cite[(2.4), Lemma~2.3]{CaoLiWu} give
	\begin{equation*}
		\|K_j^\Psi\|_1\lesssim_d\|\Psi\|_1,
		\qquad
		|\widehat K_j^\Psi(\xi)|\lesssim_d\|\Psi\|_\infty
		\min\{1,(2^j|\xi|)^{-1}\}
		\quad(\Psi\in L_\infty).
	\end{equation*}
	For $s\ge n_0$, let $(E_\nu^s)_\nu$ and $e_\nu^s\in E_\nu^s$
	be the angular partition and centers in \cite[(2.11)]{CaoLiWu}:
	\[
	\operatorname{diam}(E_\nu^s)\le2^{-s},\quad
	|e_\nu^s-e_{\nu'}^s|\ge c_d2^{-s}\ (\nu\ne\nu'),\quad
	\#\{\nu\}\lesssim_d2^{s(d-1)}.
	\]
	Take an even $\psi\in C_c^\infty(\R)$ with $0\le\psi\le1$,
	$\psi=1$ on $[-2,2]$ and support in $[-4,4]$. Set
	\begin{equation}\label{eq:Seegerkernels}
		\begin{gathered}
			\widehat{P_\nu^s h}(\xi)
			=\psi\left(2^{s/2}\frac{\langle e_\nu^s,\xi\rangle}{|\xi|}\right)
			\widehat h(\xi),\\
			H_{j\nu}^{s,\Psi}(x)=K_j^\Psi(x)\one_{E_\nu^s}(x/|x|),
			\qquad \Gamma_j^{s,\Psi}=\sum_\nu P_\nu^sH_{j\nu}^{s,\Psi}.
		\end{gathered}
	\end{equation}
	\begin{lemma}[\cite{Seeger,CaoLiWu}]\label{lem:angular-overlap}
		At each nonzero frequency, the displayed multiplier symbols obey
		\begin{equation}\label{eq:cap-fourier-overlap}
			\sup_{\xi\ne0}\sum_\nu
			\left|\psi\left(2^{s/2}
			\frac{\langle e_\nu^s,\xi\rangle}{|\xi|}\right)\right|^2
			\lesssim_d 2^{s(d-3/2)}.
		\end{equation}
		Accordingly,
		\cite[Lemma~2.5]{CaoLiWu} gives
		\[\left\|\sum_\nu P_\nu^s v_\nu\right\|_2^2
		\lesssim_d2^{s(d-3/2)}\sum_\nu\|v_\nu\|_2^2.\]
		
	\end{lemma}
	We use the remainder estimate
	\cite[Lemma~2.6]{CaoLiWu}:
	\begin{lemma}\label{lem:remainder}
		Let $s\ge n_0$, $j\in\Z$ and $\Psi\in L_\infty(S^{d-1})$.
		If $b_Q\in L_1(\A)$ is supported on a cube of side at most
		$2^{j-s+1}$ and has integral zero, then
		\[\|(K_j^\Psi-\Gamma_j^{s,\Psi})*b_Q\|_1
		\lesssim_d\|\Psi\|_\infty2^{-s/4}\|b_Q\|_1.\]
	\end{lemma}

	Fix $\delta=1/128$ and put
	\begin{equation*}
		\Omega_s^\lo=\Omega\one_{\{|\Omega|\le2^{\delta s}\|\Omega\|_1\}},
		\qquad
		\Omega_s^\hi=\Omega\one_{\{|\Omega|>2^{\delta s}\|\Omega\|_1\}}.
	\end{equation*}

	\subsection{Smooth dyadic truncations}\label{subsec:smooth-truncations}
	
	Define
	\[
	\rho(t)=\sum_{v\ge0}\beta(2^{-v}t)^2,
	\qquad S_af=\sum_{j\ge a}T_j^\Omega f.
	\]
	Then $\rho=0$ on $(0,1/2]$, $\rho=1$ on $[1,\infty)$,
	and
	\begin{equation}\label{eq:smooth-tail}
		S_af=\int_{\R^d}K_\Omega(y)\rho(2^{-a}|y|)f(\cdot-y)\dd y
		=\int_0^\infty\rho'(t)T_{\Omega,2^at}f\dd t.
	\end{equation}
The maximal theorem of Hong, Lai and Xu
	\cite[  Section~6.1]{HLX} gives
	\begin{equation*}
		\|(S_ah)_{a\in\Z}\|_{L_2(\A;\ell_\infty)}
		\lesssim_d C_{\Omega,1}\|h\|_2.
	\end{equation*}

	\section{Proof of the main theorem}\label{sec:main-proof}
	
	We prove Theorem~\ref{thm:main} first under
	\eqref{eq:finitemodel}, assuming several lemmas. Their proofs are given in the subsequent sections.
	
	\subsection{Reduction to the smooth truncations}
	
	Write
	\[
	T_{\Omega,\varepsilon}f(x)
	=\int_{|y|>\varepsilon}K_\Omega(y)f(x-y)\dd y,
	\qquad \varepsilon>0.
	\]
	For each $\varepsilon>0$, choose $a\in\Z$ with
	$2^a\le\varepsilon<2^{a+1}$. By
	\eqref{eq:dyadickernels} and \eqref{eq:smooth-tail},
	\[
	\begin{aligned}
		T_{\Omega,\varepsilon}f
		&= S_af+(T_{\Omega,\varepsilon}-S_a)f,\\
		S_af
		&=\sum_{j\ge a}K_j^\Omega*f
		=\int_{\R^d}K_\Omega(y)\rho(2^{-a}|y|)f(\cdot-y)\dd y.
	\end{aligned}
	\]
	By the support of $\beta$ and the definition of
	$\rho$, we have
	\begin{equation}\label{eq:sharp-smooth-comparison}
		\left|K_\Omega(y)
		\bigl[\one_{\{|y|>\varepsilon\}}
		-\rho(2^{-a}|y|)\bigr]\right|
		\le\sum_{j=a-1}^{a+1}|K_j^\Omega(y)|,
		\qquad 2^a\le\varepsilon<2^{a+1}.
	\end{equation}
	Consequently, it suffices to first establish the dyadic maximal estimate:
	for every $\lambda>0$,   there exists a
	projection $e\in\A$ such that
	\[
	\sup_{a\in\Z}\|eS_af e\|_\infty\le\lambda,
	\qquad
	\lambda\varphi(1-e)\le C_dC_{\Omega,1}\|f\|_1.
	\]
	
	\subsection{The   good function}
	
	The case $\Omega=0$ is immediate. Assume $\Omega\ne0$ and fix
	$\lambda>0$. Apply Theorem~\ref{thm:regularcz} to $f$ at height
	\[
	\alpha=\frac{\lambda}{C_{\Omega,1}},
	\qquad f=g+\bd+\boff.
	\]
	For every $a\in\Z$ we have the exact decomposition
	\begin{equation}\label{eq:main-direct-split}
		S_af=S_ag+S_a\bd+S_a\boff.
	\end{equation}
	Let $\zeta_0$ be the localization projection from
	Lemma~\ref{lem:cz-localization}. Then
	\[
	\varphi(1-\zeta_0)\le C_d\alpha^{-1}\|f\|_1.
	\]
	For every subprojection $e\le\zeta_0$, the identity
	\eqref{eq:main-direct-split} becomes
	\[
		eS_af e
		=eS_ag e
		+e\zeta_0 S_a\bd\zeta_0e
		+e\zeta_0 S_a\boff\zeta_0e.
	\]
	We estimate these three families simultaneously in $a$ and then
	intersect the resulting projections.
	
	The maximal
	$L_2$ estimate in Subsection~\ref{subsec:smooth-truncations} and \eqref{eq:czbounds} yield
	\[
	\begin{aligned}
		\|(S_ag)_a\|_{L_2(\A;\ell_\infty)}^2
		&\le C_dC_{\Omega,1}^2\|g\|_2^2\\
		&\le C_dC_{\Omega,1}^2\|g\|_\infty\|g\|_1
		\le C_dC_{\Omega,1}^2\alpha\|f\|_1.
	\end{aligned}
	\]
	By Lemma~\ref{lem:L2-cut}, there is a projection $e_1$ with
	\begin{equation}\label{eq:main-good-cut}
		\sup_a\|e_1S_ag e_1\|_\infty\le\lambda/8,\qquad
		\varphi(1-e_1)\le
		C_d\frac{C_{\Omega,1}^2\alpha\|f\|_1}{\lambda^2}.
	\end{equation}

	\subsection{The bad functions}
	
	We now turn to the bad functions. Put $b_{k,0}=b_{\mathrm d,k}$.
	By the decomposition in Theorem~\ref{thm:regularcz},
	\[
	\begin{aligned}
		S_a\bd &=\sum_{j\ge a}K_j^\Omega*\bd
		=\sum_{j\ge a}\sum_{k\in\Z}K_j^\Omega*b_{\mathrm d,k},\\
		S_a\boff &=\sum_{j\ge a}K_j^\Omega*\boff
		=\sum_{j\ge a}\sum_{k\in\Z}\sum_{\ell\ge1}
		K_j^\Omega*(b_{k,\ell}+g_{k,\ell}).
	\end{aligned}
	\]
	By \eqref{eq:localsupport},
	writing $n=j+k$  gives
	\begin{equation}\label{eq:main-bad-reindex}
		\begin{aligned}
			\zeta_0S_a\bd\zeta_0
			&=\zeta_0\left(\sum_{n\ge n_0}\sum_{j\ge a}
			K_j^\Omega*b_{\mathrm d,n-j}\right)\zeta_0,\\
			\zeta_0S_a\boff\zeta_0
			&=\zeta_0\left(\sum_{n\ge n_0}\sum_{\ell\ge1}\sum_{j\ge a}
			K_j^\Omega*(b_{n-j,\ell}+g_{n-j,\ell})\right)\zeta_0.
		\end{aligned}
	\end{equation}

	Let $\mathscr S=\{n_0 2^v:v\ge0\}$. Every integer $n+\ell\ge n_0$
	belongs to a unique interval $[s,2s)$ with $s\in\mathscr S$.
	For each such $s$, use
	$\Omega=\Omega_s^\lo+\Omega_s^\hi$ in \eqref{eq:main-bad-reindex}.
	The large angular terms remain combined as $u_j^s$ in
	\eqref{eq:combined-block}; the small angular terms are separated
	into the $b$- and $g$-functions. Indeed, for each $j$ and $s$,
	$$
	K_j^\Omega*u_j^s
	=K_j^{\Omega_s^\hi}*u_j^s
	+K_j^{\Omega_s^\lo}*u_j^s,$$  $$
	K_j^{\Omega_s^\lo}*u_j^s
	=\sum_{s\le n<2s}K_j^{\Omega_s^\lo}*b_{\mathrm d,n-j} 
	+\sum_{\substack{n\ge n_0,\ \ell\ge1\\s\le n+\ell<2s}}
	K_j^{\Omega_s^\lo}*(b_{n-j,\ell}+g_{n-j,\ell}).
	$$
	Summing in $j\ge a$ and $s\in\mathscr S$ gives
	\begin{equation}\label{eq:main-reconstruction}
		\begin{split}
			\zeta_0S_a(f-g)\zeta_0
			=\zeta_0\Bigg[&\sum_{s\in\mathscr S}\sum_{j\ge a}
			K_j^{\Omega_s^\hi}*u_j^s\\
			&+\sum_{s\in\mathscr S}
			\sum_{\substack{n\ge n_0,\,\ell\ge0\\s\le n+\ell<2s}}
			\sum_{j\ge a}K_j^{\Omega_s^\lo}*b_{n-j,\ell}\\
			&+\sum_{s\in\mathscr S}
			\sum_{\substack{n\ge n_0,\,\ell\ge1\\s\le n+\ell<2s}}
			\sum_{j\ge a}K_j^{\Omega_s^\lo}*g_{n-j,\ell}
			\Bigg]\zeta_0.
		\end{split}
	\end{equation}
	We estimate the three lines separately.
	The following count will be used for the last two:
	\begin{equation}\label{eq:main-pair-count}
		\#\{(n,\ell):n\ge n_0,\ \ell\ge0,\ s\le n+\ell<2s\}
		=\sum_{r=s}^{2s-1}(r-n_0+1)\le2s^2.
	\end{equation}
	
	\subsubsection{The first line of \eqref{eq:main-reconstruction}.}
	
	Young's inequality and Lemma~\ref{lem:combined-input} give
	\[
	\begin{aligned}
		\sum_{s\in\mathscr S}\sum_j\|K_j^{\Omega_s^\hi}*u_j^s\|_1
		&\le C_d\sum_{s\in\mathscr S}\|\Omega_s^\hi\|_1
		\sum_j\|u_j^s\|_1\\
		&\le C_d\|f\|_1\sum_{s\in\mathscr S}s
		\int_{\{|\Omega|>2^{s/128}\|\Omega\|_1\}}
		|\Omega(\theta)|\dd\sigma(\theta).
	\end{aligned}
	\]
	The angular summation of \cite[Section~6]{CaoLiWu} reads
	\[
	\begin{aligned}
		&\sum_{s\in\mathscr S}s
		\int_{\{|\Omega|>2^{s/128}\|\Omega\|_1\}}|\Omega|\dd\sigma\\
		&\quad=\int_{S^{d-1}}|\Omega(\theta)|
		\sum_{\substack{s\in\mathscr S\\
				2^{s/128}\|\Omega\|_1<|\Omega(\theta)|}}s\dd\sigma(\theta)
		\le C C_{\Omega,1}.
	\end{aligned}
	\]
	Apply Lemma~\ref{lem:absolute-L1} to the family
	$(K_j^{\Omega_s^\hi}*u_j^s)_{s,j}$. There is a projection $e_2$ with
	\begin{equation}\label{eq:main-high-cut}
		\begin{gathered}
			\sup_a\left\|e_2\left(\sum_{s\in\mathscr S}\sum_{j\ge a}
			K_j^{\Omega_s^\hi}*u_j^s\right)e_2\right\|_\infty
			\le\lambda/8,\\
			\varphi(1-e_2)\le C_d\frac{C_{\Omega,1}\|f\|_1}{\lambda}.
		\end{gathered}
	\end{equation}
	
	\subsubsection{ The second line of
		\eqref{eq:main-reconstruction} }
	For $\ell=0$, its individual terms are
	\[
	\sum_{j\ge a}K_j^{\Omega_s^\lo}*b_{n-j,0}
	=\sum_{j\ge a}K_j^{\Omega_s^\lo}*b_{\mathrm d,n-j};
	\]
	for $\ell\ge1$, they involve $b_{n-j,\ell}$.

	We first give a maximal estimate for partial sums.
	
	Let $h_1,\ldots,h_N\in\A$ be projections with
	$\sum_i\varphi(h_i)<\infty$, and let
	$x_i\in L_1(\A)\cap L_2(\A)$ satisfy $h_ix_i=x_i$.
	For $z\in\T$, define
	\[
	U_0(z)=1,\qquad
	U_i(z)=(1-h_1+zh_1)\cdots(1-h_i+zh_i),
	\]
	and
	\[
	\mathcal F(z)=\sum_{i=1}^N U_{i-1}(z)x_i.
	\]
	
	\begin{lemma}\label{lem:ordered-products}
		For every $\lambda>0$ and $L\ge2$,
		there exists a projection $e\in\A$ such that
		\[
		\sup_{0\le n\le N}
		\left\|e\left(\sum_{i=1}^n x_i\right)e\right\|_\infty
		\le\lambda
		\]
		and
		\[
		\begin{aligned}
			\varphi(1-e)\le C\bigg[
			&\frac1L\sum_{i=1}^N\varphi(h_i)
			+\frac1{\lambda\sqrt L}
			\sum_{i=1}^N\|x_i\|_1\\
			&+\frac{(1+\log L)^2}{\lambda^2}
			\sup_{z\in\T}\|\mathcal F(z)\|_2^2
			\bigg].
		\end{aligned}
		\]
		Here $C$ is an absolute constant.
	\end{lemma}
	
	We next introduce the enlarged dyadic grids and
	the associated projections.
	For each $\gamma\in(\Z/9\Z)^d$, let
	$\mathcal Q_k^\gamma$ be the partition into
	enlarged cubes $9I$ constructed in
	Section~\ref{sec:geometry}.
	Let $\mathcal G_k^\gamma$ be the algebra of
	functions constant on each $9I$, and let
	$\chi_k^\gamma$ be the indicator of the
	union of the selected inner cubes $I$.
	Denote by $\E_k^{\mathcal G^\gamma}$ the
	conditional expectation onto $\mathcal G_k^\gamma$:
	\[
	(\E_k^{\mathcal G^\gamma}v)(x)
	=\frac1{|Q|}\int_Q v(y)\,dy,
	\qquad x\in Q\in\mathcal Q_k^\gamma.
	\]
	For $Q\in\Q_m$, let $Q^{\rm out,\gamma}$ be the union
	of cells of $\mathcal Q_{m+4}^\gamma$ meeting $2^{n_0+10}Q$.
	Define
	\[
		\Pi_m^\gamma=\bigvee_{Q\in\Q_m}\pi_Q\one_{Q^{\rm out,\gamma}},
		\qquad
		r_k^\gamma=\bigwedge_{m\le k-4}(1-\Pi_m^\gamma).
	\]
	Set
	\begin{equation}\label{eq:final-spatial-exclusion}
		\zeta=\bigwedge_{\gamma,k}r_k^\gamma.
	\end{equation}
	
	\begin{lemma}\label{lem:spatial-projections}
		We have
		\[
		\mathcal G_k^\gamma\subset\A_k,
		\qquad
		\sum_\gamma\chi_k^\gamma=1,
		\]
		and
		\[
		r_k^\gamma\in\mathcal G_k^\gamma,\qquad
		r_{k+1}^\gamma\le r_k^\gamma,\qquad
		r_k^\gamma p_m=0\quad(m\le k-4).
		\]
		Moreover,
		\[
		\zeta\le\zeta_0,\qquad
		\varphi(1-\zeta)
		\le C_d\alpha^{-1}\|f\|_1.
		\]
	\end{lemma}
	Fix $n\ge n_0$, $\ell\ge0$, and put $r=n+\ell$.
	Let $J_j$ be scalar kernels satisfying
	\[
	\supp J_j\subset\{|x|\le2^{j+1}\},
	\qquad
	\sup_j\|J_j\|_1\le C_0\|\Omega\|_1.
	\]
	For each $\gamma$ and $k$, define
	\begin{align}
		z_k^\gamma
		&=r_k^\gamma
		[J_{-k}*(\chi_k^\gamma b_{k+n,\ell})]
		r_k^\gamma,
		\label{eq:localized-output}\\
		h_k^\gamma
		&=\supp\!\left(
		\E_k^{\mathcal G^\gamma}
		[r_k^\gamma\chi_k^\gamma
		(p_{k+n}+p_{k+r})r_k^\gamma]
		\right).
		\label{eq:localized-support}
	\end{align}
	
	\begin{lemma}\label{lem:localization}
		We have
		\[
		h_k^\gamma\in\mathcal G_k^\gamma,\qquad
		h_k^\gamma\le r_k^\gamma,\qquad
		h_k^\gamma z_k^\gamma
		=z_k^\gamma h_k^\gamma=z_k^\gamma.
		\]
		Moreover,
		\[
		\begin{aligned}
			\sum_k\varphi(h_k^\gamma)
			&\le C_d2^{dr}\alpha^{-1}\|f\|_1,\\
			\sum_k\|z_k^\gamma\|_1
			&\le C_dC_0\|\Omega\|_1\|f\|_1.
		\end{aligned}
		\]
		Finally,
		\begin{equation}\label{eq:localization-recovery}
			\zeta[J_{-k}*b_{k+n,\ell}]\zeta
			=\sum_\gamma\zeta z_k^\gamma\zeta.
		\end{equation}
	\end{lemma}
	
	Fix $\gamma$ and omit it from the notation.
	Choose integers $S\ge r+4$ and $k_0$.
	List the integers $k=k_0+mS$ ($m\in\Z$) for which $z_k\ne0$
	as $k_1<\cdots<k_N$.
	For $1\le i\le N$, write
	\[
	h_i=h_{k_i},\qquad r_i=r_{k_i},\qquad
	\chi_i=\chi_{k_i},
	\qquad
	m_i=k_i+n,\qquad l_i=k_i+r.
	\]
	Use the ordered products $U_i(z)$ defined above
	for these projections $h_i$.
	
	For either $p_i=p_{m_i}$ or $p_i=p_{l_i}$,
	set
	\[
	u_i=r_i\chi_i p_i r_i.
	\]
	
	\begin{lemma}\label{lem:localized-orthogonality}
		For $i\ne j$,
		\[
		u_iu_j=0.
		\]
		Moreover, for every $z\in\T$,
		\begin{equation}\label{eq:common-terminal-source}
			U_{i-1}(z)r_i\chi_i p_i
			=z^{-1}U_N(z)r_i\chi_i p_i.
		\end{equation}
		Consequently,
		\[
		U_{i-1}(z)u_iU_{i-1}(z)^*
		=U_N(z)u_iU_N(z)^*,
		\]
		and these operators are pairwise orthogonal.
	\end{lemma}

	Retain the projection $\zeta$ of Lemma~\ref{lem:spatial-projections}.
	Fix integers $n\ge n_0$, $\ell\ge0$ and $s\ge n_0$ such that
	$s\le r=n+\ell<2s$. We need the next estimate.
	
	\begin{lemma}\label{lem:bad-max}
		For every $t>0$, there is $e\le\zeta$ such that
		\begin{equation}\label{eq:pair-max-height}
			\sup_{a\in\Z}\left\|e\left(\sum_{j\ge a}
			K_j^{\Omega_s^\lo}*b_{n-j,\ell}\right)e\right\|_\infty\le t
		\end{equation}
		and
		\begin{equation}\label{eq:pair-max-trace}
			\begin{aligned}
				\varphi(\zeta-e)\le C_d\bigg[&
				(1+r)2^{-s}\frac{\|f\|_1}{\alpha}
				+(1+r)^22^{-s}\frac{\|\Omega\|_1\|f\|_1}{t}\\
				&+(1+r)^3s^22^{-7s/64}
				\frac{\|\Omega\|_1^2\alpha\|f\|_1}{t^2}\bigg].
			\end{aligned}
		\end{equation}
	\end{lemma}
	
	For $s\in\mathscr S$, choose
	\begin{equation}\label{eq:assigned-height}
		t_s=c_d\lambda(1+s)^{-4}2^{-s/512},\qquad
		2\sum_{s\in\mathscr S}s^2t_s\le\lambda/8,
	\end{equation}
	where $c_d>0$ is fixed sufficiently small. For each
	$n\ge n_0$, $\ell\ge0$ with $s\le n+\ell<2s$, apply
	Lemma~\ref{lem:bad-max} at height $t_s$. Denote the resulting
	projection by $e_{s,n,\ell}\le\zeta$. Thus
	\begin{equation}\label{eq:main-individual-bad-cut}
		\sup_a\left\|e_{s,n,\ell}\left(\sum_{j\ge a}
		K_j^{\Omega_s^\lo}*b_{n-j,\ell}\right)e_{s,n,\ell}\right\|_\infty
		\le t_s,
	\end{equation}
	and \eqref{eq:pair-max-trace}, with $r=n+\ell$, bounds
	$\varphi(\zeta-e_{s,n,\ell})$.
	
	Set
	\[
	e_3=\zeta\wedge\bigwedge_{s\in\mathscr S}
	\bigwedge_{\substack{n\ge n_0,\,\ell\ge0\\s\le n+\ell<2s}}
	e_{s,n,\ell}.
	\]
	Compression to a subprojection preserves each bound in
	\eqref{eq:main-individual-bad-cut}. Therefore
	\begin{equation}\label{eq:main-bad-height-sum}
		\begin{aligned}
			&\sup_a\left\|e_3\left(\sum_{s\in\mathscr S}
			\sum_{\substack{n\ge n_0,\,\ell\ge0\\s\le n+\ell<2s}}
			\sum_{j\ge a}K_j^{\Omega_s^\lo}*b_{n-j,\ell}\right)e_3
			\right\|_\infty\\
			&\quad\le\sum_{s\in\mathscr S}
			\sum_{\substack{n\ge n_0,\,\ell\ge0\\s\le n+\ell<2s}}t_s
			\le2\sum_{s\in\mathscr S}s^2t_s\le\lambda/8.
		\end{aligned}
	\end{equation}
	
	For $s\le r=n+\ell<2s$, we have $1+r\le2(1+s)$.
	Summing \eqref{eq:pair-max-trace} over the pairs counted in
	\eqref{eq:main-pair-count} gives
	\[
	\begin{aligned}
		\sum_{\substack{n\ge n_0,\,\ell\ge0\\s\le n+\ell<2s}}
		\varphi(\zeta-e_{s,n,\ell}) 
		&\le C_d(1+s)^3 2^{-s}\frac{\|f\|_1}{\alpha}
		+C_d(1+s)^4 2^{-s}\frac{\|\Omega\|_1\|f\|_1}{t_s}\\
		& +C_d(1+s)^7 2^{-7s/64}
		\frac{\|\Omega\|_1^2\alpha\|f\|_1}{t_s^2}.
	\end{aligned}
	\]
	Here the factors $(1+s)^3$, $(1+s)^4$ and $(1+s)^7$ include,
	respectively, the factors $1+r$, $(1+r)^2$ and $(1+r)^3s^2$
	from \eqref{eq:pair-max-trace} and the count $2s^2$.
	Substituting \eqref{eq:assigned-height}, we use
	\[
	\begin{gathered}
		t_s^{-1}=c_d^{-1}\lambda^{-1}(1+s)^4 2^{s/512},\qquad
		t_s^{-2}=c_d^{-2}\lambda^{-2}(1+s)^8 2^{s/256},\\
		2^{-s}2^{s/512}=2^{-511s/512},\qquad
		2^{-7s/64}2^{s/256}=2^{-27s/256}.
	\end{gathered}
	\]
	Consequently,
	\begin{equation}\label{eq:main-bad-fixed-scale}
		\begin{aligned}
			\sum_{\substack{n\ge n_0,\,\ell\ge0\\s\le n+\ell<2s}}
			\varphi(\zeta-e_{s,n,\ell}) 
			&\le C_d(1+s)^3 2^{-s}\frac{\|f\|_1}{\alpha}
			+C_d(1+s)^8 2^{-511s/512}
			\frac{\|\Omega\|_1\|f\|_1}{\lambda}\\
			& +C_d(1+s)^{15}2^{-27s/256}
			\frac{\|\Omega\|_1^2\alpha\|f\|_1}{\lambda^2}.
		\end{aligned}
	\end{equation}
	Since all $e_{s,n,\ell}$ are subprojections of the same $\zeta$,
	\[
	\zeta-e_3=\bigvee_{s\in\mathscr S}
	\bigvee_{\substack{n\ge n_0,\,\ell\ge0\\s\le n+\ell<2s}}
	(\zeta-e_{s,n,\ell}).
	\]
	Trace subadditivity, Lemma~\ref{lem:spatial-projections}, and
	\eqref{eq:main-bad-fixed-scale} imply
	\[
	\begin{aligned}
		\varphi(1-e_3)
		&\le\varphi(1-\zeta)+\sum_{s\in\mathscr S}
		\sum_{\substack{n\ge n_0,\,\ell\ge0\\s\le n+\ell<2s}}
		\varphi(\zeta-e_{s,n,\ell})\\
		&\le C_d\frac{\|f\|_1}{\alpha}
		\left(1+\sum_{s\in\mathscr S}(1+s)^3 2^{-s}\right)\\
		&\quad+C_d\frac{\|\Omega\|_1\|f\|_1}{\lambda}
		\sum_{s\in\mathscr S}(1+s)^8 2^{-511s/512}\\
		&\quad+C_d\frac{\|\Omega\|_1^2\alpha\|f\|_1}{\lambda^2}
		\sum_{s\in\mathscr S}(1+s)^{15}2^{-27s/256}.
	\end{aligned}
	\]
	The three series converge. We have proved
	\begin{equation}\label{eq:main-bad-cut}
		\varphi(1-e_3)\le C_d\left(
		\frac{\|f\|_1}{\alpha}
		+\frac{\|\Omega\|_1\|f\|_1}{\lambda}
		+\frac{\|\Omega\|_1^2\alpha\|f\|_1}{\lambda^2}\right).
	\end{equation}
	
	\subsubsection{The third line of \eqref{eq:main-reconstruction}}  Write
	\[
	K_j^{\Omega_s^\lo}*g_{n-j,\ell}
	=\Gamma_j^{s,\Omega_s^\lo}*g_{n-j,\ell}
	+(K_j^{\Omega_s^\lo}-\Gamma_j^{s,\Omega_s^\lo})*g_{n-j,\ell}.
	\]
	We estimate the first term in the maximal $L_2$ norm and the second
	by an absolutely convergent $L_1$ sum.
	For fixed $n,\ell$, the martingale difference property in
	\eqref{eq:lag-block-properties} yields the following maximal $L_2$ estimate.
	
	\begin{lemma}\label{lem:off-good-max}
		For $\Psi\in L_\infty(S^{d-1})$, $s,n\ge n_0$ and $\ell\ge1$,
		\begin{equation}\label{eq:off-good-max}
			\left\|\left(\sum_{j\ge a}\Gamma_j^{s,\Psi}*g_{n-j,\ell}
			\right)_{a\in\Z}\right\|_{L_2(\A;\ell_\infty)}
			\le C_d\|\Psi\|_\infty2^{-s/2}(\alpha\|f\|_1)^{1/2}.
		\end{equation}
	\end{lemma}
	By Lemma~\ref{lem:off-good-max} and
	$\|\Omega_s^\lo\|_\infty\le2^{s/128}\|\Omega\|_1$,
	\[
	\left\|\left(\sum_{j\ge a}
	\Gamma_j^{s,\Omega_s^\lo}*g_{n-j,\ell}\right)_a
	\right\|_{L_2(\A;\ell_\infty)}
	\le C_d\|\Omega\|_1 2^{-63s/128}(\alpha\|f\|_1)^{1/2}.
	\]
	Using \eqref{eq:main-pair-count} and the triangle inequality,
	\begin{equation}\label{eq:main-good-microlocal-sum}
		\begin{aligned}
			&\left\|\left(\sum_{s\in\mathscr S}
			\sum_{\substack{n\ge n_0,\,\ell\ge1\\s\le n+\ell<2s}}
			\sum_{j\ge a}\Gamma_j^{s,\Omega_s^\lo}*g_{n-j,\ell}\right)_a
			\right\|_{L_2(\A;\ell_\infty)}\\
			&\quad\le C_d\|\Omega\|_1(\alpha\|f\|_1)^{1/2}
			\sum_{s\in\mathscr S}2s^2 2^{-63s/128}
			\le C_d\|\Omega\|_1(\alpha\|f\|_1)^{1/2}.
		\end{aligned}
	\end{equation}
	For the second term, put $r=n+\ell$. By \eqref{eq:lag-block-properties},
	\[
	g_{n-j,\ell}=\sum_{Q\in\Q_{r-j-1}}\one_Qg_{n-j,\ell},\qquad
	\int_Qg_{n-j,\ell}=0,\qquad
	\ell(Q)=2^{j-r+1}\le2^{j-s+1}.
	\]
	Lemma~\ref{lem:remainder} therefore yields
	\[
	\begin{aligned}
		&\sum_j\|(K_j^{\Omega_s^\lo}-\Gamma_j^{s,\Omega_s^\lo})*
		g_{n-j,\ell}\|_1\\
		&\quad\le C_d\|\Omega_s^\lo\|_\infty2^{-s/4}
		\sum_j\sum_{Q\in\Q_{r-j-1}}\|\one_Qg_{n-j,\ell}\|_1\\
		&\quad=C_d\|\Omega_s^\lo\|_\infty2^{-s/4}
		\sum_j\|g_{n-j,\ell}\|_1
		\le C_d\|\Omega\|_1 2^{-31s/128}\|f\|_1.
	\end{aligned}
	\]
	Another application of \eqref{eq:main-pair-count} gives
	\begin{equation}\label{eq:main-good-remainder-sum}
		\begin{aligned}
			&\sum_{s\in\mathscr S}
			\sum_{\substack{n\ge n_0,\,\ell\ge1\\s\le n+\ell<2s}}
			\sum_j\|(K_j^{\Omega_s^\lo}-\Gamma_j^{s,\Omega_s^\lo})*
			g_{n-j,\ell}\|_1\\
			&\quad\le C_d\|\Omega\|_1\|f\|_1
			\sum_{s\in\mathscr S}2s^2 2^{-31s/128}
			\le C_d\|\Omega\|_1\|f\|_1.
		\end{aligned}
	\end{equation}
	Apply Lemma~\ref{lem:L2-cut} to
	\eqref{eq:main-good-microlocal-sum} and
	Lemma~\ref{lem:absolute-L1} to \eqref{eq:main-good-remainder-sum}.
	We obtain projections $e_4,e_5$ with
	\begin{equation}\label{eq:main-off-good-cuts}
		\begin{gathered}
			\sup_a\left\|e_4\left(\sum_{s\in\mathscr S}
			\sum_{\substack{n\ge n_0,\,\ell\ge1\\s\le n+\ell<2s}}
			\sum_{j\ge a}\Gamma_j^{s,\Omega_s^\lo}*g_{n-j,\ell}
			\right)e_4\right\|_\infty\le\lambda/8,\\
			\sup_a\left\|e_5\left(\sum_{s\in\mathscr S}
			\sum_{\substack{n\ge n_0,\,\ell\ge1\\s\le n+\ell<2s}}
			\sum_{j\ge a}(K_j^{\Omega_s^\lo}-\Gamma_j^{s,\Omega_s^\lo})*
			g_{n-j,\ell}\right)e_5\right\|_\infty\le\lambda/8,\\
			\varphi(1-e_4)\le C_d\frac{\|\Omega\|_1^2\alpha\|f\|_1}{\lambda^2},
			\qquad
			\varphi(1-e_5)\le C_d\frac{\|\Omega\|_1\|f\|_1}{\lambda}.
		\end{gathered}
	\end{equation}
	
	\subsection{Completion of the bad functions}
	
	Put $e=e_1\wedge\cdots\wedge e_5$. Since $e\le e_3\le\zeta\le\zeta_0$,
	\eqref{eq:main-reconstruction} remains valid after compression by $e$.
	Equations \eqref{eq:main-good-cut}, \eqref{eq:main-high-cut},
	\eqref{eq:main-off-good-cuts} and \eqref{eq:main-bad-height-sum} give
	\[
	\sup_a\|eS_afe\|_\infty\le5\lambda/8\le\lambda.
	\]
	Combining their trace estimates with \eqref{eq:main-bad-cut},
	\[
	\begin{aligned}
		\varphi(1-e)
		\le\sum_{i=1}^5\varphi(1-e_i) 
		&\le C_d\left[
		\frac1\alpha+\frac{C_{\Omega,1}+\|\Omega\|_1}{\lambda}
		+\frac{(C_{\Omega,1}^2+\|\Omega\|_1^2)\alpha}{\lambda^2}
		\right]\|f\|_1\\
		&\le C_d\frac{C_{\Omega,1}\|f\|_1}{\lambda},
	\end{aligned}
	\]
	where $\alpha=\lambda/C_{\Omega,1}$ and $\|\Omega\|_1\le C_{\Omega,1}$.
	Thus
	\begin{equation}\label{eq:smooth-endpoint}
		\sup_{a\in\Z}\|eS_afe\|_\infty\le\lambda,\qquad
		\lambda\varphi(1-e)\le C_dC_{\Omega,1}\|f\|_1.
	\end{equation}

	\subsection{Completion of Theorem~\ref{thm:main}}
	
	We first consider $f\ge0$ satisfying
	\eqref{eq:finitemodel}. For $K=K_\Omega$,
	Seeger's conditions give
	\[
	V=(\log2)|\Omega|,
	\qquad
	\eta(a)\lesssim_d a^{-1}\|\Omega\|_1.
	\]
	By Lemma~\ref{lem:positive-annuli}, for every
	$t>0$ there is a projection $e$ such that
	\[
	\sup_j\|e(|K_j^\Omega|*f)e\|_\infty\le t,
	\qquad
	\varphi(1-e)\le
	C_d\frac{C_{\Omega,1}\|f\|_1}{t}.
	\]
	For $2^a\le\varepsilon<2^{a+1}$,
	\eqref{eq:sharp-smooth-comparison} gives
	\[
	\|e(T_{\Omega,\varepsilon}-S_a)fe\|_\infty
	\le2\sum_{j=a-1}^{a+1}
	\|e(|K_j^\Omega|*f)e\|_\infty.
	\]
	Intersecting the projections from
	\eqref{eq:smooth-endpoint} and
	Lemma~\ref{lem:positive-annuli}, at heights
	$\lambda/2$ and $\lambda/12$, respectively,
	we obtain
	\[
	\sup_{\varepsilon>0}
	\|eT_{\Omega,\varepsilon}fe\|_\infty
	\le\lambda,
	\qquad
	\varphi(1-e)\le
	C_d\frac{C_{\Omega,1}\|f\|_1}{\lambda}.
	\]
	Decomposing into four positive parts extends
	this estimate to the linear span of
	\eqref{eq:finitemodel}.
	The case of general $f\in L_1(\A)$ follows by density argument.
	This completes the proof of
	Theorem~\ref{thm:main}.

	\Needspace{10\baselineskip}
	\section{Proof of Lemma~\ref{lem:ordered-products}}\label{sec:products}
	
	\begin{proof}
		Fix $\lambda>0$ and $L\ge2$. Set
		\[
		\eta=1-L^{-2},\qquad a=b=\lambda/2,
		\qquad W_i(z)=1-h_i+zh_i.
		\]
		Let $m$ be normalized Haar measure on $\T$ and set
		\[
		c_\eta=\int_\T\frac{dm(z)}{|1-\eta\bar z|}.
		\]
		Put $\mathcal F_n(z)=\sum_{i\le n}U_{i-1}(z)x_i$ and
		$C_\eta(z)=(1-\eta\bar z)^{-1}$. The factors $W_i(z)$ are
		unitary for $z\in\T$ and contractive for $z=\eta$. Since $h_ix_i=x_i$,
		\[
		U_n(z)^*U_{i-1}(z)x_i=
		\begin{cases}
			\bar z\,W_n(z)^*\cdots W_{i+1}(z)^*x_i,&i\le n,\\
			W_{n+1}(z)\cdots W_{i-1}(z)x_i,&i>n.
		\end{cases}
		\]
		The first expression contains only negative
		powers of $z$, while the second is a polynomial
		in $z$. Thus Cauchy's formula gives
		\[
		\mathcal F_n(\eta)
		=\int_\T C_\eta(z)
		[1-U_n(\eta)U_n(z)^*]\mathcal F(z)\,dm(z).
		\]
		Set
		\[
		B_n(z)=1-U_n(\eta)U_n(z)^*.
		\]
		Since $\|B_n(z)\|_\infty\le2$,
		by Jensen's inequality,
		\[
		\begin{aligned}
			|\mathcal F_n(\eta)|^2
			&=\left|
			\int_\T C_\eta(z)B_n(z)\mathcal F(z)\,dm(z)
			\right|^2\\
			&\le c_\eta\int_\T |C_\eta(z)|
			\mathcal F(z)^*B_n(z)^*B_n(z)
			\mathcal F(z)\,dm(z)\\
			&\le4c_\eta\int_\T|C_\eta(z)|
			\mathcal F(z)^*\mathcal F(z)\,dm(z)
			=:Q_\eta.
		\end{aligned}
		\]
		Taking traces, we obtain
		\[
		\begin{aligned}
			\varphi(Q_\eta)
			&=4c_\eta\int_\T|C_\eta(z)|
			\|\mathcal F(z)\|_2^2\,dm(z)\\
			&\le4c_\eta^2
			\sup_{z\in\T}\|\mathcal F(z)\|_2^2.
		\end{aligned}
		\]
		Thus $e_2=\one_{[0,b^2]}(Q_\eta)$ satisfies
		\[
			\sup_n\|\mathcal F_n(\eta)e_2\|_\infty\le b,
			\qquad
			\varphi(1-e_2)\le\frac{4c_\eta^2}{b^2}
			\sup_{z\in\T}\|\mathcal F(z)\|_2^2.
		\]
				Set
		$V_{i,n}=W_{i+1}\cdots W_n$. In the following estimates, write
		$W_i=W_i(\eta)$ and $U_i=U_i(\eta)$. The identities
		$W_i^2=1-(1-\eta^2)h_i$ give
		$$
		1-U_n =(1-\eta)\sum_{i\le n}U_{i-1}h_i
		=(1-\eta)\sum_{i\le n}h_iV_{i,n}, $$ $$
		\sum_{i\le n}U_{i-1}h_iU_{i-1}^*
		=\frac{1-U_nU_n^*}{1-\eta^2}\le\frac1{1-\eta^2}, $$ $$
		\sum_{i\le n}V_{i,n}^*h_iV_{i,n}
		=\frac{1-U_n^*U_n}{1-\eta^2}\le\frac1{1-\eta^2}. $$
		Applying the discrete  case of \cite[Proposition 1.1]{Mei}   to the two identities
		for $1-U_n$
		gives, for every projection $q$,
		\begin{equation*}
			\max\{\|(1-U_n)q\|_\infty^2,\|q(1-U_n)\|_\infty^2\}
			\le\frac{1-\eta}{1+\eta}
			\left\|q\left(\sum_{i\le n}h_i\right)q\right\|_\infty.
		\end{equation*}
		Take $e_0=\one_{[0,L]}(\sum_i h_i)$ and
		$v_i=e_0[1-U_{i-1}(\eta)]x_ie_0$. Then
		\[
		\varphi(1-e_0)\le L^{-1}\sum_i\varphi(h_i),\qquad
		\sum_i\|v_i\|_1
		\le\left(\frac{1-\eta}{1+\eta}L\right)^{1/2}\sum_i\|x_i\|_1.
		\]
		Lemma~\ref{lem:absolute-L1}, applied in $e_0\A e_0$, provides
		$e_1\le e_0$ such that
		$$
		\sup_n\left\|e_1\left(\sum_{i\le n}v_i\right)e_1\right\|_\infty
		\le a,$$ $$
		\varphi(e_0-e_1)
		\le\frac2a\left(\frac{1-\eta}{1+\eta}L\right)^{1/2}
		\sum_i\|x_i\|_1.
		$$
		For $e=e_1\wedge e_2$,
		\[
		e\left(\sum_{i\le n}x_i\right)e
		=e\mathcal F_n(\eta)e+e\left(\sum_{i\le n}v_i\right)e,
		\quad
		\varphi(1-e)\le\varphi(1-e_0)+\varphi(e_0-e_1)+\varphi(1-e_2).
		\]
		Since $a=b=\lambda/2$, the height bound is $\lambda$. Moreover,
		\[
		\left(\frac{1-\eta}{1+\eta}L\right)^{1/2}
		\le L^{-1/2},\qquad c_\eta\lesssim1+\log L.
		\]
		The trace estimate follows from the three preceding bounds,
		with an absolute constant independent of $N$.
	\end{proof}
	
	\Needspace{10\baselineskip}
	\section{Proofs of Lemmas \ref{lem:spatial-projections}, \ref{lem:localization},  \ref{lem:localized-orthogonality}}\label{sec:geometry}

	Let $\mathbf{1}_d=(1,\ldots,1)$.
	For each $\gamma\in(\Z/9\Z)^d$, define
	\[
	a_0^\gamma=\gamma,\qquad
	a_{j+1}^\gamma=5a_j^\gamma+2\mathbf{1}_d,
	\quad j\in\Z,
	\]
	where all operations are taken modulo $9$.
	
	At scale $j$, select the cubes $I=2^j(v+[0,1)^d)$ with
	$v\equiv a_j^\gamma\pmod9$. Their enlargements $9I$ partition
	$\R^d$. The congruence
	$2(a_{j+1}^\gamma-4\mathbf{1}_d)\equiv a_j^\gamma-4\mathbf{1}_d\pmod9$
	makes these partitions nested. With $k=-j$, denote the resulting
	partition by $\mathcal Q_k^\gamma$, its algebra by
	$\mathcal G_k^\gamma$, and the indicator of the selected inner
	cubes by $\chi_k^\gamma$. Then
	\[
	\mathcal G_k^\gamma\subset\A_k,\qquad
	\sum_\gamma\chi_k^\gamma=1.
	\]
	If $\supp J_{-k}\subset\{|x|\le2^{-k+1}\}$, then
	$\supp[J_{-k}*(\one_Iv)]\subset9I$. Hence, for
	$A,B\in\mathcal G_k^\gamma$,
	\begin{equation}\label{eq:masks-through-convolution}
		A[J_{-k}*(\chi_k^\gamma v)]B
		=\sum_{I\text{ selected}}J_{-k}*(\one_IAvB)
		=J_{-k}*(\chi_k^\gamma AvB).
	\end{equation}
	
	\begin{proof}[Proof of Lemma~\ref{lem:spatial-projections}]
		The grid construction gives
		$\mathcal G_k^\gamma\subset\A_k$ and
		$\sum_\gamma\chi_k^\gamma=1$.
		Since $|Q^{\rm out,\gamma}|\le C_d|Q|$ and
		$Q\subset2^{n_0+10}Q\subset Q^{\rm out,\gamma}$,
		\[
		r_k^\gamma\in\mathcal G_k^\gamma,\quad
		r_{k+1}^\gamma\le r_k^\gamma,\quad
		r_k^\gamma p_m=0\ (m\le k-4),\quad r_k^\gamma\le q_{k-4}.
		\]
		Thus $\zeta\le\zeta_0$, and
		\begin{equation*}
			\varphi(1-\zeta)
			\le\sum_{\gamma,m,Q\in\Q_m}|Q^{\rm out,\gamma}|\tau(\pi_Q)
			\le C_d\sum_m\varphi(p_m)
			\le C_d\alpha^{-1}\|f\|_1.
		\end{equation*}
	\end{proof}

	\begin{proof}[Proof of Lemma~\ref{lem:localization}]
		Fix a grid and suppress $\gamma$. Write
		$a_k=r_k\chi_k(p_{k+n}+p_{k+r})r_k$. The definition of $h_k$
		and bimodularity give $h_k\le r_k$ and
		\[
		\varphi((1-h_k)a_k(1-h_k))
		=\varphi((1-h_k)\E_k^{\mathcal G}(a_k)(1-h_k))=0.
		\]
		Consequently,
		\begin{equation*}
			h_kr_k\chi_kp_{k+n}=r_k\chi_kp_{k+n},\qquad
			h_kr_k\chi_kp_{k+r}=r_k\chi_kp_{k+r}.
		\end{equation*}
		Together with \eqref{eq:masks-through-convolution}, this proves
		$h_kz_k=z_kh_k=z_k$.
		
		On a cell $K=9I\in\mathcal Q_k^\gamma$,
		\[
		(\E_k^{\mathcal G}a_k)|_K
		=\frac1{|K|}\sum_{m=k+n,k+r}
		\sum_{\substack{Q\in\Q_m\\Q\subset I}}|Q|r_K\pi_Qr_K.
		\]
		By polar decomposition,
		$\tau(\supp(r_K\pi_Qr_K))\le\tau(\pi_Q)$. Thus
		\[
		\begin{aligned}
			\sum_k\varphi(h_k)
			&\le9^d(2^{dn}+2^{dr})\sum_m\varphi(p_m)
			\le C_d2^{dr}\alpha^{-1}\|f\|_1,\\
			\sum_k\|z_k\|_1
			&\le C_0\|\Omega\|_1\sum_k\|b_{k+n,\ell}\|_1
			\le C_dC_0\|\Omega\|_1\|f\|_1.
		\end{aligned}
		\]
		When $\ell=0$, the identical support terms are counted twice,
		which only changes the constant. Finally,
		\[
		\zeta[J_{-k}*b_{k+n,\ell}]\zeta
		=\sum_\gamma\zeta r_k^\gamma
		[J_{-k}*(\chi_k^\gamma b_{k+n,\ell})]r_k^\gamma\zeta
		=\sum_\gamma\zeta z_k^\gamma\zeta.
		\]
	\end{proof}
	
	\begin{proof}[Proof of Lemma~\ref{lem:localized-orthogonality}]
		Write $W_i(z)=1-h_i+zh_i$.
		For $i<j$ in one residue class, $k_j\ge k_i+S\ge l_i+4$.
		Therefore
		\begin{equation*}
			r_j\le r_i,\qquad r_jp_{m_i}=r_jp_{l_i}=0.
		\end{equation*}
		For either choice $p_i=p_{m_i}$ or $p_i=p_{l_i}$,
		\[
		u_iu_j=\chi_i\chi_jr_ip_ir_ir_jp_jr_j
		=\chi_i\chi_jr_ip_ir_jp_jr_j=0,
		\qquad \supp(u_i)\supp(u_j)=0.
		\]
		Moreover, for $j>i$,
		\[
		h_{k_j}r_i\chi_ip_i=\chi_i h_{k_j}p_i=0,
		\qquad W_i(z)r_i\chi_ip_i=z r_i\chi_ip_i.
		\]
		Reading the product from the right gives
		\[
		U_N(z)r_i\chi_ip_i
		=U_{i-1}(z)W_i(z)r_i\chi_ip_i
		=zU_{i-1}(z)r_i\chi_ip_i.
		\]
		This is \eqref{eq:common-terminal-source}. In particular,
		\[
		U_{i-1}(z)u_iU_{i-1}(z)^*=U_N(z)u_iU_N(z)^*,
		\]
		so the conjugated range projections remain pairwise orthogonal.
	\end{proof}

	\Needspace{10\baselineskip}
	\section{Proof of Lemma~\ref{lem:bad-max}}\label{sec:microlocal}
In this section, we prove Lemma~\ref{lem:bad-max}. We begin with some preliminary results.
 
	Use the notation of Lemma~\ref{lem:localized-orthogonality}. Fix $z\in\T$, put $s_i=\chi_i r_i$, and, for either
	$\mathcal F_i=f$ or $\mathcal F_i=f_{l_i}$, define
	\[
	\begin{aligned}
		A_i&=U_{i-1}(z)s_i p_{m_i}\mathcal F_i p_{m_i}s_iU_{i-1}(z)^*,\\
		D_i&=s_i p_{l_i}\mathcal F_i p_{l_i}s_i,\\
		X_i&=U_{i-1}(z)s_i p_{m_i}\mathcal F_i p_{l_i}s_i,\qquad
		d_i(z)=U_{i-1}(z)s_i b_{m_i,\ell}s_i.
	\end{aligned}
	\]
	\begin{lemma}\label{lem:source-corners}
		There are projections $P_i\in\A_{m_i}$ and $Q_i\in\A_{l_i}$ such that
		\[
		\begin{gathered}
			\begin{pmatrix}A_i&X_i\\X_i^*&D_i\end{pmatrix}\ge0,\qquad
			A_i=P_iA_iP_i,\quad D_i=Q_iD_iQ_i,\\
			\sum_iP_i\le1,\quad\sum_iQ_i\le1,\quad
			\E_{m_i}A_i\le2^d\alpha P_i,\quad\E_{l_i}D_i\le2^d\alpha Q_i,\\
			\sum_i\|A_i\|_1\le\|f\|_1,\qquad
			\sum_i\|D_i\|_1\le\|f\|_1.
		\end{gathered}
		\]
		The same estimates hold when $p_{m_i}$
		and $p_{l_i}$ are interchanged.
		Moreover, $\sum_i\|d_i(z)\|_1\le4\|f\|_1$ and
		\begin{equation}\label{eq:localized-product-sum}
			\sum_iU_{i-1}(z)z_{k_i}=\sum_iJ_{-k_i}*d_i(z),\qquad
			\E_{l_i}d_i(z)=0.
		\end{equation}
		For $\ell=0$, the $L_1$ bound is $2\|f\|_1$.
	\end{lemma}
	\begin{proof}
		Set
		\[
		P_i^0=\supp(s_ip_{m_i}s_i),\quad Q_i^0=\supp(s_ip_{l_i}s_i),
		\quad P_i=U_{i-1}(z)P_i^0U_{i-1}(z)^*,\quad Q_i=Q_i^0.
		\]
		Since $U_{i-1}(z),r_i\in\mathcal G_{k_i}\subset\A_{k_i}$
		and $\chi_i,s_i\in\A_{k_i}$,
		\[
		P_i\in\A_{m_i},\quad Q_i\in\A_{l_i},\quad
		P_i=U_N(z)P_i^0U_N(z)^*,\quad
		\sum_iP_i\le1,\quad\sum_iQ_i\le1.
		\]
		The matrix factorization
		\[
		\begin{pmatrix}A_i&X_i\\X_i^*&D_i\end{pmatrix}
		=\begin{pmatrix}U_{i-1}s_ip_{m_i}\\s_ip_{l_i}\end{pmatrix}
		\mathcal F_i
		\begin{pmatrix}p_{m_i}s_iU_{i-1}^*&p_{l_i}s_i\end{pmatrix}\ge0
		\]
		proves positivity and the support assertions. For
		$\mathcal F_i=f$ and for $\mathcal F_i=f_{l_i}$, bimodularity gives
		\[
		\begin{aligned}
			\E_{m_i}A_i
			&=U_{i-1}s_ip_{m_i}f_{m_i}p_{m_i}s_iU_{i-1}^*
			\le2^d\alpha P_i,\\
			\E_{l_i}D_i&=s_ip_{l_i}f_{l_i}p_{l_i}s_i
			\le2^d\alpha Q_i,\\
			\sum_i\varphi(A_i)
			&=\sum_i\varphi(fp_{m_i}s_ip_{m_i})
			\le\sum_i\varphi(fp_{m_i})\le\|f\|_1,\\
			\sum_i\varphi(D_i)
			&=\sum_i\varphi(fp_{l_i}s_ip_{l_i})\le\|f\|_1.
		\end{aligned}
		\]
		All factors multiplying $f-f_{l_i}$ belong to $\A_{l_i}$, so
		\[
		\E_{l_i}d_i(z)=U_{i-1}(z)s_i\E_{l_i}(b_{m_i,\ell})s_i=0.
		\]
		For each positive matrix, H\"older's inequality and the positive
		matrix factorization  give
		\[
		\sum_i\|X_i\|_1
		\le\sum_i\|A_i\|_1^{1/2}\|D_i\|_1^{1/2}
		\le\left(\sum_i\|A_i\|_1\right)^{1/2}
		\left(\sum_i\|D_i\|_1\right)^{1/2}\le\|f\|_1.
		\]
		The two choices of $\mathcal F_i$ and the two orientations
		therefore give $\sum_i\|d_i(z)\|_1\le4\|f\|_1$; for $\ell=0$
		there are only two terms. Finally, \eqref{eq:masks-through-convolution}
		yields
		\[
		U_{i-1}(z)z_{k_i}
		=J_{-k_i}*[U_{i-1}(z)\chi_i r_i b_{m_i,\ell}r_i]
		=J_{-k_i}*d_i(z).
		\]
	\end{proof}
	
	For a positive operator matrix
	$\left(\begin{smallmatrix}A&X\\X^*&D\end{smallmatrix}\right)$,
	the standard factorization $X=A^{1/2}vD^{1/2}$, $\|v\|_\infty\le1$,
	gives
	\begin{equation}\label{eq:block-holder}
		\|X\|_p\le\|A\|_p^{1/2}\|D\|_p^{1/2},\qquad 1\le p<\infty.
	\end{equation}
 For two such matrices,
the  Cauchy--Schwarz inequality gives
	\begin{equation}\label{eq:block-cross-trace}
		|\tau(X_1^*X_2)|
		\le[\tau(A_1A_2)\tau(D_1D_2)]^{1/2}
		\le\tfrac12\tau(A_1A_2+D_1D_2).
	\end{equation}
	
	The next lemma extends the form of the microlocal estimate in
	\cite[Lemmas~3.1 and~3.4]{CaoLiWu} to positive matrices.
	
	Fix $n\ge n_0$, $\ell\ge0$, $r=n+\ell$, and $n_0\le s\le r$.
	Consider finite families $A_j,D_j\in L_1(\A)_+\cap\A$ and $X_j\in\A$
	with
	\begin{equation}\label{eq:abstract-corner-block}
		\begin{pmatrix}A_j&X_j\\X_j^*&D_j\end{pmatrix}\ge0.
	\end{equation}
	Assume that projections $P_j\in\A_{n-j}$ and $Q_j\in\A_{r-j}$ satisfy
	\begin{equation}\label{eq:abstract-corner-averages}
		\sum_jP_j\le1,\quad\sum_jQ_j\le1,\quad
		\E_{n-j}A_j\le\alpha P_j,\quad \E_{r-j}D_j\le\alpha Q_j,
	\end{equation}
	and
	\begin{equation}\label{eq:abstract-corner-mass}
		\sum_j\|A_j\|_1\le\|f\|_1,\qquad
		\sum_j\|D_j\|_1\le\|f\|_1.
	\end{equation}
	Let $H_j$ be scalar kernels with
	\[
	\supp H_j\subset\{2^{j-2}\le|x|\le2^{j+2}\},\qquad
	|H_j(x)|\le M_0 2^{-jd}.
	\]
	Using \eqref{eq:Seegerkernels}, set
	\[
	H_{j\nu}^s(x)=H_j(x)\one_{E_\nu^s}(x/|x|),\qquad
	\Gamma_j^s=\sum_\nu P_\nu^sH_{j\nu}^s.
	\]
	For $H_j=K_j^\Psi$, this is $\Gamma_j^{s,\Psi}$.
	
	\begin{lemma}\label{lem:corner-microlocal}
		Under \eqref{eq:abstract-corner-block}--\eqref{eq:abstract-corner-mass},
		\[
		\left\|\sum_j\Gamma_j^s*X_j\right\|_2
		\lesssim_d M_0\sqrt{1+(s-n)_+}\,
		2^{-s/4}(\alpha\|f\|_1)^{1/2}.
		\]
		The same bound holds with the two conditional expectation levels
		in \eqref{eq:abstract-corner-averages} interchanged.
	\end{lemma}
	
	\begin{proof}[Proof of Lemma~\ref{lem:corner-microlocal}]
		Fix $\nu$ and put $Z_j=H_{j\nu}^s*X_j$. The scalar sector estimates in
		\cite[proof of Lemma~2.1]{Seeger} give symmetric rectangles $R_{j\nu}$,
		with side lengths comparable to $2^j$ and $2^{j-s}$, such that
		\begin{align*}
			|H_{j\nu}^s|&\le C_dM_02^{-jd}\one_{R_{j\nu}},
			&|R_{j\nu}|&\le C_d2^{jd-s(d-1)},\\
			\|H_{j\nu}^s\|_1&\le C_dM_02^{-s(d-1)},
			&|\widetilde H_{i\nu}^s*H_{j\nu}^s|
			&\le C_dM_0^22^{-jd-s(d-1)}\one_{R_{j\nu}}
			\quad(i\le j).
		\end{align*}
		Here $\widetilde H(x)=\overline{H(-x)}$; fixed dimensional dilations of
		$R_{j\nu}$ are absorbed in the constants. Set $a=(s-n)_+$. Then
		\begin{equation}\label{eq:corner-near-far}
			\left\|\sum_jZ_j\right\|_2^2
			\le(2a+1)\sum_j\|Z_j\|_2^2
			+2\sum_j\sum_{i<j-a}|\langle Z_i,Z_j\rangle_2|.
		\end{equation}
		
		The positivity of \eqref{eq:abstract-corner-block}, integrated against
		$|H_{j\nu}^s|$ with its scalar phase in the off-diagonal entries, implies
		\begin{equation}\label{eq:corner-near-block}
			\|Z_j(x)\|_{L_2(\M)}^2
			\le \tau\big((|H_{j\nu}^s|*A_j)(x)\big)
			\big\|(|H_{j\nu}^s|*D_j)(x)\big\|_\infty.
		\end{equation}
		Since $2^{j-r}\le2^{j-s}$, completing the level-$(r-j)$ cubes gives
		\[
		\begin{aligned}
			(|H_{j\nu}^s|*D_j)(x)
			&\le C_dM_02^{-jd}\int_{x+R_{j\nu}}D_j(y)\dd y\\
			&\le C_dM_02^{-jd}\int_{x+C_dR_{j\nu}}\E_{r-j}D_j(y)\dd y
			\le C_d\alpha M_02^{-s(d-1)}\one.
		\end{aligned}
		\]
		Consequently,
		\begin{equation}\label{eq:corner-diagonal-energy}
			\sum_j\|Z_j\|_2^2
			\le C_d\alpha M_02^{-s(d-1)}
			\sum_j\|H_{j\nu}^s\|_1\|A_j\|_1
			\le C_d\alpha M_0^22^{-2s(d-1)}\|f\|_1.
		\end{equation}
		
		For $i<j-a$, both source cube sizes satisfy
		$2^{i-r}\le2^{i-n}\le2^{j-s-1}$. Hence
		\[
		\begin{aligned}
			\sum_{i<j-a}\int_{x+R_{j\nu}}A_i(y)\dd y
			&\le\int_{x+C_dR_{j\nu}}\sum_{i<j-a}\E_{n-i}A_i(y)\dd y\\
			&\le\alpha\int_{x+C_dR_{j\nu}}\sum_iP_i(y)\dd y
			\le C_d\alpha2^{jd-s(d-1)}\one,\\
			\sum_{i<j-a}\int_{x+R_{j\nu}}D_i(y)\dd y
			&\le\alpha\int_{x+C_dR_{j\nu}}\sum_iQ_i(y)\dd y
			\le C_d\alpha2^{jd-s(d-1)}\one.
		\end{aligned}
		\]
		The convolution pairing and \eqref{eq:block-cross-trace} now yield
		\[
		\begin{aligned}
			\sum_{i<j-a}|\langle Z_i,Z_j\rangle_2|
			&\le\frac12\sum_{i<j-a}\iint
			|(\widetilde H_{i\nu}^s*H_{j\nu}^s)(y-x)|\\[-2mm]
			&\hspace{28mm}\times
			\tau\big(A_i(y)A_j(x)+D_i(y)D_j(x)\big)\dd y\dd x\\
			&\le C_d\alpha M_0^22^{-2s(d-1)}
			\big(\|A_j\|_1+\|D_j\|_1\big).
		\end{aligned}
		\]
		Together with \eqref{eq:corner-near-far}--\eqref{eq:corner-diagonal-energy},
		this proves
		\[
		\left\|\sum_jH_{j\nu}^s*X_j\right\|_2^2
		\le C_d(1+a)\alpha M_0^22^{-2s(d-1)}\|f\|_1.
		\]
		Lemma~\ref{lem:angular-overlap} and $\#\{\nu\}\lesssim_d2^{s(d-1)}$ give
		\[
		\begin{aligned}
			\left\|\sum_j\Gamma_j^s*X_j\right\|_2^2
			&\lesssim_d2^{s(d-3/2)}
			\sum_\nu\left\|\sum_jH_{j\nu}^s*X_j\right\|_2^2\\
			&\lesssim_d(1+a)\alpha M_0^2\|f\|_1
			2^{s(d-3/2)+s(d-1)-2s(d-1)}\\
			&= (1+a)\alpha M_0^2\|f\|_1\,2^{-s/2}.
		\end{aligned}
		\]
		Interchanging $A_j,D_j$ in \eqref{eq:corner-near-block} proves the
		assertion with the two source levels reversed.
	\end{proof}

	Retain \eqref{eq:abstract-corner-block}--\eqref{eq:abstract-corner-mass}
	and let $K_j(x,y)$ be scalar kernels satisfying
	\[
	|K_j(x,y)|\le M_02^{-jd}\one_{\{|x-y|\le C_d2^j\}}.
	\]
	\begin{lemma}\label{lem:positive-kernel-lp}
		For $1<p<\infty$,
		\[
		\left\|\sum_j\int K_j(\cdot,y)X_j(y)\dd y\right\|_p
		\lesssim_{d,p}M_0\alpha^{1-1/p}\|f\|_1^{1/p}.
		\]
		The constant is independent of $n,\ell,s$ and the number of indices.
	\end{lemma}
	
	\begin{proof}
		Set $a_j=\E_{n-j}A_j$. By the hypotheses,
		\[
		0\le\sum_ja_j\le\alpha\one,\qquad
		\left\|\sum_ja_j\right\|_1\le\|f\|_1,\qquad
		\left\|\sum_ja_j\right\|_p
		\le\alpha^{1-1/p}\|f\|_1^{1/p}.
		\]
		Completing the source cubes and applying the positive dual
		Hardy--Littlewood inequality \eqref{eq:dual-HL}, we obtain
		\[
		\begin{aligned}
			\left\|\sum_j\int |K_j(\cdot,y)|A_j(y)\dd y\right\|_p
			&\le C_dM_0\left\|\sum_j A_{C_d2^j}a_j\right\|_p\\
			&\le C_{d,p}M_0\left\|\sum_ja_j\right\|_p
			\le C_{d,p}M_0\alpha^{1-1/p}\|f\|_1^{1/p}.
		\end{aligned}
		\]
		The same bound holds for $D_j$, using $\E_{r-j}D_j$. Integration of the
		positive matrices \eqref{eq:abstract-corner-block}, with the scalar
		phase of $K_j$ in their off-diagonal entries, and
		\eqref{eq:block-holder} therefore give
		\[
		\begin{aligned}
			\left\|\sum_j\int K_j(\cdot,y)X_j(y)\dd y\right\|_p
			&\le\left\|\sum_j\int |K_j(\cdot,y)|A_j(y)\dd y\right\|_p^{1/2}
			\left\|\sum_j\int |K_j(\cdot,y)|D_j(y)\dd y\right\|_p^{1/2}\\
			&\le C_{d,p}M_0\alpha^{1-1/p}\|f\|_1^{1/p}.
		\end{aligned}
		\]
	\end{proof}
	
	In particular, for $K_j(x,y)=H_j(x-y)$,
	\[\left\|\sum_jH_j*X_j\right\|_3
	\lesssim_d M_0\alpha^{2/3}\|f\|_1^{1/3}.\]
	
	We will also use the following consequence of the distribution
	inequalities and integration formulas in
	\cite[Proposition~2.2, Lemma~2.5(v), Corollary~2.8]{FackKosaki}:
	if $Y=G+R$, with $G\in L_2(\A)$, $R\in L_1(\A)$ and
	$Y\in L_3(\A)$, then
	\begin{equation}\label{eq:mixed-interpolation}
		\|Y\|_2\lesssim \|G\|_2+
		\|R\|_1^{1/4}\|Y\|_3^{3/4}.
	\end{equation}

	Fix $n\ge n_0$, $\ell\ge0$ and $n_0\le s\le r=n+\ell$.
	Use Lemmas~\ref{lem:localization} and
	\ref{lem:localized-orthogonality} with $J_j=K_j^{\Omega_s^\lo}$. Fix a grid and one residue class
	modulo $S\ge r+4$, and list its nonzero levels as
	$k_1<\cdots<k_N$. Let $z_{k_i}$ be the localized outputs in
	\eqref{eq:localized-output}, let $h_{k_i}$ be the projections
	defined in \eqref{eq:localized-support}, and form $U_i(z)$ as in
	Lemma~\ref{lem:ordered-products}.
	
	\begin{lemma}\label{lem:transformed-full-energy}
		With this notation,
		\begin{equation*}
			\sup_{|z|=1}\left\|\sum_{i=1}^NU_{i-1}(z)z_{k_i}\right\|_2
			\le C_d \|\Omega\|_12^{-7s/128}(\alpha \|f\|_1)^{1/2}.
		\end{equation*}
	\end{lemma}
	\begin{proof}
		Fix $z\in\T$. By \eqref{eq:localized-product-sum} and
		Lemma~\ref{lem:source-corners},
		\[
		Y(z):=\sum_iU_{i-1}(z)z_{k_i}
		=\sum_iK_{-k_i}^{\Omega_s^\lo}*d_i(z),\qquad
		\E_{k_i+r}d_i(z)=0,\qquad \sum_i\|d_i(z)\|_1\le4\|f\|_1.
		\]
		Each $d_i(z)$ is a sum of at most four signed off-diagonal entries
		of the positive matrices in Lemma~\ref{lem:source-corners}.
		Apply Lemmas~\ref{lem:corner-microlocal} and~\ref{lem:positive-kernel-lp}
		to these entries. With
		$G(z)=\sum_i\Gamma_{-k_i}^{s,\Omega_s^\lo}*d_i(z)$, this gives
		\begin{align*}
			\|G(z)\|_2
			&\le C_d\|\Omega_s^\lo\|_\infty
			\sqrt{1+(s-n)_+}\,2^{-s/4}(\alpha\|f\|_1)^{1/2},
			\\
			\|Y(z)\|_3
			&\le C_d\|\Omega_s^\lo\|_\infty\alpha^{2/3}\|f\|_1^{1/3}.
		\end{align*}
		On each $Q\in\Q_{k_i+r}$, the source has integral zero and
		$\ell(Q)=2^{-k_i-r}\le2^{-k_i-s}$. Lemma~\ref{lem:remainder} yields
		\begin{equation*}
			\|Y(z)-G(z)\|_1
			\le C_d\|\Omega_s^\lo\|_\infty2^{-s/4}
			\sum_i\|d_i(z)\|_1
			\le C_d\|\Omega_s^\lo\|_\infty2^{-s/4}\|f\|_1.
		\end{equation*}
		Substituting these estimates into \eqref{eq:mixed-interpolation},
		\[
		\begin{aligned}
			\|Y(z)\|_2
			&\le C_d\|\Omega_s^\lo\|_\infty
			\big[\sqrt{1+(s-n)_+}\,2^{-s/4}+2^{-s/16}\big]
			(\alpha\|f\|_1)^{1/2}\\
			&\le C_d\|\Omega\|_1\,2^{s/128-s/16}(\alpha\|f\|_1)^{1/2}
			=C_d\|\Omega\|_1\,2^{-7s/128}(\alpha\|f\|_1)^{1/2}.
		\end{aligned}
		\]
		All bounds are uniform in $z\in\T$.
	\end{proof}
	Now we are ready to prove Lemma~\ref{lem:bad-max}.
	\begin{proof}[Proof of Lemma~\ref{lem:bad-max}]
		Use $\zeta$ from \eqref{eq:final-spatial-exclusion}, set $r=n+\ell$,
		and separate $k=-j$ into residue classes modulo $r+5$.
		For each grid and residue class, Lemma~\ref{lem:localization} gives
		\[
		\sum_i\varphi(h_{k_i})\le C_d2^{dr}\|f\|_1/\alpha,
		\qquad
		\sum_i\|z_{k_i}\|_1\le C_d\|\Omega\|_1\|f\|_1.
		\]
		Apply Lemma~\ref{lem:ordered-products} at height $u$ with
		\[
		L=2^{dr+s}.
		\]
		Since $s\le r<2s$ and $d\ge2$,
		\[
		\frac{\sum_i\varphi(h_{k_i})}{L}
		\le C_d2^{-s}\frac{\|f\|_1}{\alpha},\qquad
		L^{-1/2}\le2^{-s},\qquad
		1+\log L\lesssim_d s.
		\]
		Lemma~\ref{lem:transformed-full-energy} supplies the remaining term.
		Thus one projection $e_\iota$ controls all initial sums at height $u$ and
		\begin{equation*}
			\varphi(1-e_\iota)\le C_d\left[
			2^{-s}\frac{\|f\|_1}{\alpha}
			+2^{-s}\frac{\|\Omega\|_1\|f\|_1}{u}
			+s^22^{-7s/64}\frac{\|\Omega\|_1^2\alpha\|f\|_1}{u^2}\right].
		\end{equation*}
		The condition $j\ge a$ is $k_i\le-a$, so these are precisely the required
		tails. There are $9^d(r+5)$ choices of $\iota$. Put
		\[
		u=\frac{t}{9^d(r+5)},\qquad e=\zeta\wedge\bigwedge_\iota e_\iota.
		\]
		By \eqref{eq:localization-recovery} and trace subadditivity for joins,
		 $$
	 			\sup_a\left\|e\sum_{j\ge a}K_j^{\Omega_s^\lo}*b_{n-j,\ell}\,e\right\|_\infty
		 \le9^d(r+5)u=t$$    \begin{align*}
		 	 	\varphi(\zeta-e)
		 	  \le\sum_\iota\varphi(1-e_\iota) 
		 	 &\le C_d\left[(1+r)2^{-s}\frac{\|f\|_1}{\alpha}
		 	 +(1+r)^22^{-s}\frac{\|\Omega\|_1\|f\|_1}{t}\right]\\
		 	 &\quad+C_d(1+r)^3s^22^{-7s/64}
		 	 \frac{\|\Omega\|_1^2\alpha\|f\|_1}{t^2}.
		 \end{align*}
		  		This is \eqref{eq:pair-max-height}--\eqref{eq:pair-max-trace}.
	\end{proof}
	
	\Needspace{10\baselineskip}
	\section{Proof of Lemma~\ref{lem:off-good-max}}\label{sec:off-good}

	In this section, we prove Lemma~\ref{lem:off-good-max}. We begin with some preliminary results.
	
	Fix a real radial
	$\vartheta\in C_c^\infty(\R^d\setminus\{0\})$, supported in
	$\{1/2\le|\xi|\le2\}$, with
	$\sum_{v\in\Z}\vartheta(2^{-v}\xi)^2=1$ for $\xi\ne0$, and set
	\[
	\widehat{\Delta_vu}(\xi)
	=\vartheta(2^{-v}\xi)^2\widehat u(\xi).
	\]
	Then $\sum_v\Delta_vu=u$ in $L_2(\A)$.
	
	\begin{lemma}
		\label{lem:martingale-fourier}
		If $u\in L_2(\A)$, $\E_m u=u$ and $\E_{m-1}u=0$, then
		\begin{equation}\label{eq:MD-LP}
			\norm[2]{\Delta_vu}
			\le C_d2^{-9|m-v|/20}\norm[2]{u}.
		\end{equation}
	\end{lemma}
	
	This is a  basic operator-valued extension of  
	\cite[Proposition~4.4]{GrafakosKalton}. We omit  the proof.

	Let $0<c<C<\infty$ and let $(v_l)_{l\in\Z}\subset L_2(\A)$ satisfy
	\[
	\supp\widehat v_l\subset\{\xi:c2^l\le|\xi|\le C2^l\},\qquad
	\sum_l\|v_l\|_2^2<\infty.
	\]
	\begin{lemma}\label{lem:annular-synthesis}
		The partial sums converge in $L_2(\A)$ and satisfy
		\[
		\left\|\left(\sum_{l\le r}v_l\right)_{r\in\Z}\right\|_{L_2(\A;\ell_\infty)}
		\le C_{d,c,C}\left(\sum_l\|v_l\|_2^2\right)^{1/2}.
		\]
		The same bound holds for   $\sum_{l\ge r}v_l$.
	\end{lemma}
	\begin{proof}
		We use \cite[Proposition~2.3]{HWWFourier} at $p=2$:
		\begin{equation*}
			\|(z_r)_r\|_{L_2(\A;\ell_\infty)}
			\le\left(\sum_r\|z_r\|_2^2\right)^{1/2}.
		\end{equation*}
		Let $\Phi_r$ have symbol $\eta(2^{-r}\xi)$, where $\eta$ is a fixed
		smooth radial cutoff equal to one on the unit ball and supported in
		the ball of radius two. Its Schwartz kernel is dominated by a summable
		combination of ball averages. The noncommutative Hardy--Littlewood
		inequality \cite[Theorem~4.1]{HongLiaoWang} gives
		\begin{equation*}
			\|(\Phi_rh)_r\|_{L_2(\A;\ell_\infty)}\le C_d\|h\|_2.
		\end{equation*}
		For a finite family put $v=\sum_lv_l$. Finite overlap of the Fourier
		supports implies $\|v\|_2^2\le C_{c,C}\sum_l\|v_l\|_2^2$.
		Choose $b_0$ with $C2^{-b_0}\le1$ and $c2^{b_0}\ge2$. Then
		\[
		\sum_{l\le r}v_l=\Phi_rv+E_r,\qquad
		E_r=\sum_{|l-r|<b_0}(\one_{\{l\le r\}}\Id-\Phi_r)v_l,
		\qquad \sum_r\|E_r\|_2^2\le C_{c,C}\sum_l\|v_l\|_2^2.
		\]
		Combining the preceding estimates,
		\[
		\left\|\left(\sum_{l\le r}v_l\right)_r\right\|_{L_2(\A;\ell_\infty)}
		\le C_d\|v\|_2+\left(\sum_r\|E_r\|_2^2\right)^{1/2}
		\le C_{d,c,C}\left(\sum_l\|v_l\|_2^2\right)^{1/2}.
		\]
	 
	\end{proof}
	
Now we give the proof of Lemma~\ref{lem:off-good-max}.
	
	\begin{proof}[Proof of Lemma~\ref{lem:off-good-max}]
		The sector kernels satisfy
		\[
		\norm[L_1(\R^d)]{H_{j\nu}^{s,\Psi}}
		\le C_d\norm[\infty]{\Psi}\,2^{-s(d-1)}.
		\]
	By \eqref{eq:cap-fourier-overlap},
		\begin{align*}
			|\widehat{\Gamma_j^{s,\Psi}}(\xi)|
			&\le\sum_{\nu:\widehat P_\nu^s(\xi)\ne0}
			|\widehat P_\nu^s(\xi)|\,
			|\widehat H_{j\nu}^{s,\Psi}(\xi)|\\
			&\le C_d2^{s(d-3/2)}
			\|\Psi\|_\infty2^{-s(d-1)}
			=C_d\|\Psi\|_\infty2^{-s/2}.
		\end{align*}
		In particular,
		\begin{equation}\label{eq:Gamma-multiplier}
			\sup_{j\in\Z}\operatorname*{ess\,sup}_{\xi\in\R^d}
			|\widehat{\Gamma_j^{s,\Psi}}(\xi)|
			\le C_d\|\Psi\|_\infty2^{-s/2}.
		\end{equation}
				Put $m_j=n+\ell-j$. By \eqref{eq:lag-block-properties},
		$g_{n-j,\ell}$ is a martingale difference at level $m_j$.
		For $q\in\Z$, set
		\[
		V_j^{(q)}=\Gamma_j^{s,\Psi}*
		\Delta_{m_j+q}g_{n-j,\ell}.
		\]
		Its Fourier support lies in
		$\{2^{m_j+q-1}\le|\xi|\le2^{m_j+q+1}\}$.
		For fixed $q$, the index $l=m_j+q$ identifies $j\ge a$ with
		$l\le n+\ell+q-a$. Hence Lemma~\ref{lem:annular-synthesis},
		Plancherel's theorem, \eqref{eq:Gamma-multiplier}, and
		\eqref{eq:MD-LP} give
		\[
		\begin{split}
			\norm[L_2(\A;\ell_\infty)]{
				\left(\sum_{j\ge a}V_j^{(q)}\right)_a}
			&\le C_d\left(\sum_j\norm[2]{V_j^{(q)}}^2\right)^{1/2}\\
			&\le C_d\norm[\infty]{\Psi}\,2^{-s/2}2^{-9|q|/20}
			\left(\sum_j\norm[2]{g_{n-j,\ell}}^2\right)^{1/2}.
		\end{split}
		\]
		Since $\sum_q2^{-9|q|/20}<\infty$, summation in $q$ converges
		in $L_2(\A;\ell_\infty)$. The Fourier resolution identifies
		the component at $a$ with
		$\sum_{j\ge a}\Gamma_j^{s,\Psi}*g_{n-j,\ell}$.
		The square-sum bound in \eqref{eq:lag-block-properties} proves
		\eqref{eq:off-good-max}.
			\end{proof}
	
	\appendix
	\section{Nonhomogeneous kernels}
	\label{sec:nonhomogeneous}
	
	We prove the maximal analogue of \cite[Theorem~6.1]{CaoLiWu}
	under the nonhomogeneous kernel assumptions of
	\cite[(1.4)--(1.7)]{Seeger}. Let
	$K\in\mathcal S'(\R^d)$, $d\ge2$, be locally integrable away
	from the origin and suppose that $\widehat K\in L_\infty(\R^d)$.
	For $\theta\in S^{d-1}$ and $a\ge2$, set
	\begin{align*}
		V_R(\theta)&=\int_R^{2R}|K(r\theta)|r^{d-1}\dd r, \quad
		 V(\theta)=\sup_{R>0}V_R(\theta),
		\\
		\eta(a)&=\sup_{\substack{h>0\\R\ge ah}}
		\int_{S^{d-1}}\int_R^{2R}
		|K((r-h)\theta)-K(r\theta)|r^{d-1}\dd r\dd\sigma(\theta).\end{align*}
	Write
	\[\begin{aligned}
		C_{V,1}&=\int_{S^{d-1}}V(\theta)
		\left(1+\log^+\frac{V(\theta)}{\|V\|_1}\right)
		\dd\sigma(\theta),\\
		A_K&=\|\widehat K\|_\infty+C_{V,1}
		+\int_2^\infty\eta(a)\frac{\dd a}{a}.
	\end{aligned}\]
	As before, $C_{V,1}=0$ when $V=0$.  
	For $\varepsilon>0$, 
	define
	\[
	T_{K,\varepsilon}f(x)
	=\int_{|y|>\varepsilon}K(y)f(x-y)\dd y.
	\]
	\begin{theorem}\label{thm:general-seeger}
		Suppose that $A_K<\infty$. For every $f\in L_1(\A)$ and
		$\lambda>0$ there is a projection $e\in\A$ such that
		\begin{equation*}
			\sup_{\varepsilon>0}
			\|eT_{K,\varepsilon}f e\|_\infty\le\lambda,
			\qquad
			\lambda\varphi(1-e)\le C_dA_K\|f\|_1,
		\end{equation*}

	\end{theorem}
  Assume $V\ne0$.
		Use the partition $\beta$ from \eqref{eq:dyadickernels} and put
	\[
	K_j(x)=\beta(2^{-j}|x|)^2K(x),\qquad
	S_a^Kf=\sum_{j\ge a}K_j*f.
	\]
	Thus
	$S_a^Kf=[K(\cdot)\rho(2^{-a}|\cdot|)]*f$, with the smooth
	cutoff $\rho$ in \eqref{eq:smooth-tail}. Polar coordinates give
	\begin{equation}\label{eq:general-annular-L1}
		\sup_j\|K_j\|_1\le C_d\|V\|_1.
	\end{equation}

	Set $\delta_0=[100(d+2)]^{-1}$. Choose a nonnegative
	$\vartheta_0\in C_c^\infty((0,1/16))$ with integral one.
	For every integer $t\ge1$, use the radial regularization from
	\cite[Section~6]{CaoLiWu}:
	 $$
		J_j^t(r\theta)
		 =\beta(2^{-j}r)^2\int_0^\infty 
		K((r-u)\theta)2^{\delta_0t-j}
		\vartheta_0(2^{\delta_0t-j}u)\dd u,$$  $$
		H_j^t(r\theta)
	 =\one_{\{V\le2^{\delta_0t}\|V\|_1\}}(\theta)J_j^t(r\theta), $$
	$$	R_j^t=J_j^t-H_j^t,\qquad S_j^t=K_j-J_j^t.
$$
	\begin{lemma}\label{lem:general-smoothing}
		The kernels $H_j^t,J_j^t,R_j^t,S_j^t$ are supported in
		$\{2^{j-1}<|x|<2^{j+1}\}$. For every integer $b\ge0$,
		\begin{equation}\label{eq:general-derivatives}
			\sup_{j,r,\theta}r^{d+b}
			|\partial_r^bH_j^t(r\theta)|
			\le C_{d,b}\|V\|_1\,2^{\delta_0t(b+2)}.
		\end{equation}
		Moreover,
		\begin{align}
			&\sup_j\bigl(\|J_j^t\|_1+\|H_j^t\|_1\bigr)\le C_d\|V\|_1,
			\label{eq:general-smoothing-L1}\\
				&\sup_j\|R_j^t\|_1
		\le C_d\int_{\{V>2^{\delta_0t}\|V\|_1\}}V\dd\sigma,\notag
			\\
			&\sup_j\|S_j^t\|_1
			\le C_d\eta(2^{\delta_0t+2}).\notag
		\end{align}
	\end{lemma}
	The radial regularization estimates \cite[(6.5)--(6.7)]{CaoLiWu}
	hold for every $t\ge1$; the dyadic restriction there is used in
	the subsequent summation. The additional bound
	\eqref{eq:general-smoothing-L1} follows from Fubini's theorem
	and the defining radial shell bound for $V$.

	For $q\ge1$, define
	\[
	E_q=C_d\left[
	\|V\|_1\,2^{-(1-4\delta_0)q}
	+\int_{\{V>2^{\delta_0q}\|V\|_1\}}V\dd\sigma
	+\eta(2^{\delta_0q+2})\right].
	\]
	\begin{lemma}\label{lem:general-fourier-envelope}
		The constant in $E_q$ can be chosen so that
		\[
		\begin{gathered}
			|\widehat K_j(\xi)|\le E_q
			\quad(2^{q-1}\le2^j|\xi|\le2^{q+1}),\\
			\sum_{q\ge1}E_q\le C_d\left[C_{V,1}
			+\int_2^\infty\eta(a)\frac{\dd a}{a}\right].
		\end{gathered}
		\]
		Both bounds also hold with $K$ replaced by $|K|$.
	\end{lemma}
	\begin{proof}
		Apply Lemma~\ref{lem:general-smoothing} with $t=q$. The scalar
		Fourier estimate in \cite[proof of Lemma~2.3]{CaoLiWu}, applied with
		the radial derivative bounds \eqref{eq:general-derivatives}, gives
		\[
		|\widehat H_j^q(\xi)|
		\le C_d\|V\|_1\,2^{4\delta_0q}(2^j|\xi|)^{-1}
		\le C_d\|V\|_1\,2^{-(1-4\delta_0)q}.
		\]
		Therefore
		\[
		|\widehat K_j(\xi)|
		\le|\widehat H_j^q(\xi)|+\|R_j^q\|_1+\|S_j^q\|_1\le E_q.
		\]
		The angular counting and Dini summation in
		\cite[Section~6]{CaoLiWu} yield
		\[
		\sum_{q\ge1}E_q
		\le C_d\left[C_{V,1}+\int_2^\infty\eta(a)\frac{\dd a}{a}\right].
		\]
		Since $V_{|K|}=V$ and $\eta_{|K|}\le\eta_K$, the same estimates hold
		for $|K|$.
	\end{proof}
	
	\subsection{The good function}
	
	The boundedness of $\widehat K$ supplies a bound for the full
	convolution operator. The following argument supplies the
	maximal estimate required for the good function.
	
	Let $(L_j)_{j\in\Z}$ be integrable scalar kernels and let
	$(d_q)_{q\in\Z}$ be a summable nonnegative sequence satisfying
	\[
	|\widehat L_j(\xi)|\le d_q
	\quad(2^{q-1}\le2^j|\xi|\le2^{q+1}).
	\]
	\begin{lemma}\label{lem:multiplier-envelope-max}
		For $h\in L_2(\A)$, the tails converge in $L_2(\A)$ and
		\[
		\left\|\left(\sum_{j\ge a}L_j*h\right)_{a\in\Z}
		\right\|_{L_2(\A;\ell_\infty)}
		\le C_d\left(\sum_qd_q\right)\|h\|_2.
		\]
	\end{lemma}
	\begin{proof}
		Use $\Delta_v$ from Section~\ref{sec:off-good}. For fixed $q$,
		$L_j*\Delta_{q-j}h$ has Fourier support at scale $2^{q-j}$ and
		\[
		\sum_j\|L_j*\Delta_{q-j}h\|_2^2
		\le d_q^2\sum_j\|\Delta_{q-j}h\|_2^2
		\le C_dd_q^2\|h\|_2^2.
		\]
		Lemma~\ref{lem:annular-synthesis}, with $l=q-j$, gives
		\[
		\begin{aligned}
			\left\|\left(\sum_{j\ge a}L_j*h\right)_a\right\|_{L_2(\A;\ell_\infty)}
			&\le\sum_q\left\|\left(\sum_{j\ge a}L_j*\Delta_{q-j}h\right)_a
			\right\|_{L_2(\A;\ell_\infty)}\\
			&\le C_d\sum_qd_q\|h\|_2.
		\end{aligned}
		\]
	\end{proof}
	
	Fix a nonnegative radial
	$\phi\in C_c^\infty(\{1/2<|x|<2\})$ with integral one, and
	write
	\[
	\phi_j(x)=2^{-jd}\phi(2^{-j}x).
	\]
	These normalized annular averages are bounded by a constant
	times the positive ball averages at radius $2^{j+1}$.
	
	\begin{lemma}\label{lem:general-good-max}
		For every $h\in L_2(\A)$,
		\begin{equation*}
			\|(S_a^Kh)_{a\in\Z}\|_{L_2(\A;\ell_\infty)}
			\le C_dA_K\|h\|_2.
		\end{equation*}
	\end{lemma}
	\begin{proof}
		Put $\kappa_j=\int K_j$ and $L_j=K_j-\kappa_j\phi_j$.
		The support, cancellation and \eqref{eq:general-annular-L1} give
		$|\widehat L_j(\xi)|\lesssim_d\|V\|_1\,2^j|\xi|$ when
		$2^j|\xi|\le1$. At relative frequency $2^q$, $q\ge1$,
		Lemma~\ref{lem:general-fourier-envelope} gives
		\[
		|\widehat L_j(\xi)|\le E_q+C_d\|V\|_1\,2^{-2q}.
		\]
		Thus Lemma~\ref{lem:multiplier-envelope-max} implies
		\begin{equation}\label{eq:general-cancelled-tails}
			\left\|\left(\sum_{j\ge a}L_j*h\right)_a\right\|_{L_2(\A;\ell_\infty)}
			\le C_d\left[C_{V,1}+\int_2^\infty\eta(a)\frac{\dd a}{a}\right]\|h\|_2.
		\end{equation}
		
		For the annular means, take $\chi(x)=1-\rho(|x|)\in C_c^\infty(\R^d)$.
		Fourier duality and dilation invariance of the Fourier $L_1$ norm give
		\begin{equation*}
			\left|\sum_{j=a}^b\kappa_j\right|
			=\left|\left\langle K,
			\chi(2^{-b-1}\cdot)-\chi(2^{-a}\cdot)\right\rangle\right|
			\le C_d\|\widehat K\|_\infty.
		\end{equation*}
		Choose $c_j-c_{j-1}=\kappa_j$ with
		$\sup_j|c_j|\le C_d\|\widehat K\|_\infty$. Abel summation gives
		\[
		\sum_{j=a}^b\kappa_j(\phi_j*h)
		=c_b(\phi_b*h)-c_{a-1}(\phi_a*h)
		+\sum_{j=a}^{b-1}c_j(\phi_j-\phi_{j+1})*h.
		\]
		The first term tends to zero in $L_2$ as $b\to\infty$.
		The noncommutative Hardy--Littlewood inequality controls
		$(c_{a-1}\phi_a*h)_a$. Moreover,
		\[
		|c_j|\,|\widehat\phi(2^j\xi)-\widehat\phi(2^{j+1}\xi)|
		\le C_d\|\widehat K\|_\infty
		\min\{2^j|\xi|,(2^j|\xi|)^{-2}\}.
		\]
		A second application of Lemma~\ref{lem:multiplier-envelope-max} therefore
		controls the remaining tails by $C_d\|\widehat K\|_\infty\|h\|_2$.
		Together with \eqref{eq:general-cancelled-tails}, this proves the lemma.
	\end{proof}

	The passage to sharp radii uses a common projection estimate for
	the positive annular kernels $|K_j|$.
	
	\begin{lemma}\label{lem:positive-annuli}
		Let $K$ be locally integrable away from the origin with
		$A_K<\infty$. For every
		$f\in L_1(\A)_+$ and $t>0$, there is a projection $e$ such that
		\begin{equation}\label{eq:positive-annuli}
			\sup_{j\in\Z}\|e(|K_j|*f)e\|_\infty\le t,
			\qquad
			\varphi(1-e)\le
			\frac{C_d\|f\|_1}{t}
			\left(C_{V,1}+\int_2^\infty\eta(a)\frac{\dd a}{a}\right).
		\end{equation}
	\end{lemma}
	\begin{proof}
		Put $\mu_j=\int|K_j|$ and $N_j=|K_j|-\mu_j\phi_j$.
		Then $0\le\mu_j\le C_d\|V\|_1$ and $\int N_j=0$.
		Lemma~\ref{lem:general-fourier-envelope}, applied to $|K|$, gives
		\begin{equation*}
			\sup_{\xi\ne0}\sum_j|\widehat N_j(\xi)|
			\le C_d\left[C_{V,1}+\int_2^\infty\eta(a)\frac{\dd a}{a}\right].
		\end{equation*}
		For a finite set $J\subset\Z$ and signs $\epsilon_j\in\{-1,1\}$,
		let $N_\epsilon=\sum_{j\in J}\epsilon_jN_j$. The scalar cutoff
		$\sum_{j\in J}\epsilon_j\beta(2^{-j}r)^2$ has modulus at most one
		and derivative bounded by $C/r$. Hence, uniformly in $J,\epsilon$,
		\[
		V_{N_\epsilon}(\theta)\le C_d(V(\theta)+\|V\|_1),\qquad
		\eta_{N_\epsilon}(a)\le C_d\eta(a)+C_d\|V\|_1/a,
		\qquad C_{V_{N_\epsilon},1}\le C_dC_{V,1}.
		\]
		For the last estimate one uses the scalar inequality
		\[
		0\le U\le W\quad\Longrightarrow\quad
		\int U\log^+\frac{U}{\|U\|_1}
		\le\int W\log^+\frac{W}{\|W\|_1}+\|W\|_1/e.
		\]
		The general-kernel theorem \cite[Theorem~6.1]{CaoLiWu} now yields
		\[
		\|N_\epsilon*f\|_{1,\infty}
		\le C_d\left[C_{V,1}+\int_2^\infty\eta(a)\frac{\dd a}{a}\right]\|f\|_1.
		\]
		For independent Rademacher signs the same bound holds in the tensor
		product algebra. The weak noncommutative Khintchine inequality
		\cite[Corollary~3.2]{Cadilhac} supplies $N_j*f=u_j+v_j$ with
		\[
		\left\|\left(\sum_{j\in J}u_j^*u_j\right)^{1/2}\right\|_{1,\infty}
		+\left\|\left(\sum_{j\in J}v_jv_j^*\right)^{1/2}\right\|_{1,\infty}
		\le C_d\left[C_{V,1}+\int_2^\infty\eta(a)\frac{\dd a}{a}\right]\|f\|_1.
		\]
		Take the intersection $e$ of the spectral projections at height $t$
		of these two square functions. Then
		\[
		\begin{aligned}
			\sup_{j\in J}\|e(N_j*f)e\|_\infty
			&\le\sup_{j\in J}(\|u_je\|_\infty+\|ev_j\|_\infty)\le2t,\\
			\varphi(1-e)
			&\le\frac{C_d\|f\|_1}{t}
			\left[C_{V,1}+\int_2^\infty\eta(a)\frac{\dd a}{a}\right].
		\end{aligned}
		\]
		Since $|K_j|*f=N_j*f+\mu_j\phi_j*f$ and
		$0\le\mu_j\phi_j*f\le C_d\|V\|_1 A_{2^{j+1}}f$,
		the Hardy--Littlewood weak $(1,1)$ inequality
		\cite[Theorem~4.1]{HongLiaoWang} gives the same bound for the finite
		positive family $(|K_j|*f)_{j\in J}$. 
		The passage to \eqref{eq:positive-annuli} for all $j\in\Z$ follows by standard arguments.
	\end{proof}
	
	\subsection{The bad function}
	
	We record the radial derivative loss in the microlocal
	remainder before using the ordered products.
	
	For $s\ge n_0$, form $\Gamma_j^s$ from $H_j^s$ using
	\eqref{eq:Seegerkernels}.
	\begin{lemma}\label{lem:general-microlocal-remainder}
		If $b_Q\in L_1(\A)$ has integral zero and is supported in a cube
		of side at most $2^{j-s+1}$, then
		\begin{equation*}
			\|(H_j^s-\Gamma_j^s)*b_Q\|_1
			\le C_d\|V\|_1\,2^{-s/8}\|b_Q\|_1.
		\end{equation*}
	\end{lemma}
	\begin{proof}
		Use \cite[Lemma~2.2]{Seeger}, with $\kappa=1/2$,
		$\varepsilon=1/4$ and $N=4d+4$. Substituting
		\eqref{eq:general-derivatives} into the cited estimate gives
		\[
		\begin{aligned}
			\|(H_j^s-\Gamma_j^s)*b_Q\|_1
			&\le C_d\|V\|_1
			\left[2^{(-1/4+2\delta_0)s}
			+2^{[-1+(4d+6)\delta_0]s}\right]\|b_Q\|_1\\
			&\le C_d\|V\|_1\,2^{-s/8}\|b_Q\|_1,
			\qquad \delta_0=[100(d+2)]^{-1}.
		\end{aligned}
		\]
	\end{proof}
	
	Let $f$ satisfy \eqref{eq:finitemodel}, use its decomposition at
	height $\alpha$, and retain $\zeta$ from Lemma~\ref{lem:spatial-projections}.
	Fix $n\ge n_0$, $\ell\ge0$ and $n_0\le s\le r=n+\ell<2s$.
	\begin{lemma}\label{lem:general-bad-max}
		For every $u>0$, there is $e\le\zeta$ such that
		\begin{equation}\label{eq:general-pair-height}
			\sup_a\left\|e\left(\sum_{j\ge a}
			H_j^s*b_{n-j,\ell}\right)e\right\|_\infty\le u
		\end{equation}
		and
		\begin{equation}\label{eq:general-pair-trace}
			\begin{split}
				\varphi(\zeta-e)\le C_d\bigg[
				&(1+r)2^{-s}\frac{\|f\|_1}{\alpha}
				+(1+r)^22^{-s}\frac{\|V\|_1\|f\|_1}{u}\\
				&+(1+r)^3s^22^{-s/32}\frac{\|V\|_1^2\alpha \|f\|_1}{u^2}\bigg].
			\end{split}
		\end{equation}
	\end{lemma}
	\begin{proof}
		Apply Lemmas~\ref{lem:localization} and
		\ref{lem:localized-orthogonality} with $J_j=H_j^s$.
		For a fixed grid and residue class use the sources of
		Lemma~\ref{lem:source-corners}, so that
		\[
		Y(z)=\sum_iU_{i-1}(z)z_{k_i}=\sum_iH_{-k_i}^s*d_i(z),\qquad
		G(z)=\sum_i\Gamma_{-k_i}^s*d_i(z).
		\]
		Lemmas~\ref{lem:corner-microlocal}, \ref{lem:positive-kernel-lp}
		and~\ref{lem:general-microlocal-remainder} give
		\[
		\begin{aligned}
			\|G(z)\|_2
			&\le C_d\|V\|_1\,2^{2\delta_0s-s/4}
			\sqrt{1+(s-n)_+}(\alpha\|f\|_1)^{1/2},\\
			\|Y(z)\|_3
			&\le C_d\|V\|_1\,2^{2\delta_0s}\alpha^{2/3}\|f\|_1^{1/3},\\
			\|Y(z)-G(z)\|_1
			&\le C_d\|V\|_1\,2^{-s/8}\|f\|_1.
		\end{aligned}
		\]
		By \eqref{eq:mixed-interpolation},
		\begin{equation*}
			\begin{aligned}
				\sup_{z\in\T}\|Y(z)\|_2
				&\le C_d\|V\|_1(\alpha\|f\|_1)^{1/2}
				\left[2^{(2\delta_0-1/4)s}\sqrt{1+(s-n)_+}
				+2^{(-1/32+3\delta_0/2)s}\right]\\
				&\le C_d\|V\|_1\,2^{-s/64}(\alpha\|f\|_1)^{1/2}.
			\end{aligned}
		\end{equation*}
		The proof of Lemma~\ref{lem:bad-max} now applies with
		$\|\Omega\|_1$ replaced by $\|V\|_1$ and the squared decay
		$2^{-7s/64}$ replaced by $2^{-s/32}$. Specifically, apply Lemma~\ref{lem:ordered-products} with
		$L=2^{dr+s}$ and height $u/[9^d(r+5)]$
		for each grid and residue class. The support and $L_1$
		bounds are unchanged by \eqref{eq:general-smoothing-L1}; hence
		\eqref{eq:general-pair-height}--\eqref{eq:general-pair-trace} follow.
	\end{proof}
	
	\subsection{Completion of the proof}
	
	\begin{proof}[Proof of Theorem~\ref{thm:general-seeger}]
	We first assume \eqref{eq:finitemodel}. Apply the
		Calder\'on--Zygmund decomposition at height $\alpha>0$.
		In \eqref{eq:main-reconstruction}, replace $K_j^{\Omega_s^\lo}$ by
		$H_j^s$ and $K_j^{\Omega_s^\hi}$ by $R_j^s+S_j^s$.
		The combined inputs $u_j^s$ remain those in \eqref{eq:combined-block}.
		By Lemma~\ref{lem:combined-input}, the regularization estimates, and
		the scalar summations of \cite[Section~6]{CaoLiWu},
		\begin{equation*}
			\begin{aligned}
				\sum_{s\in\mathscr S}\sum_j\|(R_j^s+S_j^s)*u_j^s\|_1
				&\lesssim_d\|f\|_1\sum_{s\in\mathscr S}s
				\left[\int_{\{V>2^{\delta_0s}\|V\|_1\}}V\dd\sigma
				+\eta(2^{\delta_0s+2})\right]\\
				&\lesssim_d A_K\|f\|_1.
			\end{aligned}
		\end{equation*}
		Lemma~\ref{lem:absolute-L1} controls all tails of this sum.
		
		For the terms containing $H_j^s*b_{n-j,\ell}$, apply
		Lemma~\ref{lem:general-bad-max} with the heights $t_s$ in
		\eqref{eq:assigned-height}. The counting is the same as in
		Section~\ref{sec:main-proof}. The three trace sums are bounded by
		\[
		\begin{aligned}
			&C_d\frac{\|f\|_1}{\alpha}\sum_s(1+s)^3 2^{-s}
			+C_d\frac{\|V\|_1\|f\|_1}{\lambda}
			\sum_s(1+s)^8 2^{-511s/512}\\
			&\qquad+C_d\frac{\|V\|_1^2\alpha\|f\|_1}{\lambda^2}
			\sum_s(1+s)^{15}2^{-7s/256}\\
			&\hspace{12mm}\le C_d\left[
			\frac{\|f\|_1}{\alpha}
			+\frac{\|V\|_1\|f\|_1}{\lambda}
			+\frac{\|V\|_1^2\alpha\|f\|_1}{\lambda^2}\right].
		\end{aligned}
		\]
		Here $2^{-s/32}2^{s/256}=2^{-7s/256}$.
		
		For $g_{n-j,\ell}$, split $H_j^s=\Gamma_j^s+(H_j^s-\Gamma_j^s)$.
		The proof of Lemma~\ref{lem:off-good-max} uses only the sector multiplier
		bound; \eqref{eq:general-derivatives} therefore gives
		\[
		\left\|\left(\sum_{j\ge a}\Gamma_j^s*g_{n-j,\ell}\right)_a
		\right\|_{L_2(\A;\ell_\infty)}
		\le C_d\|V\|_1\,2^{(2\delta_0-1/2)s}(\alpha\|f\|_1)^{1/2}.
		\]
		For the complement, use the parent-cube cancellation of $g_{n-j,\ell}$
		and Lemma~\ref{lem:general-microlocal-remainder}. Summing over the at most
		$2s^2$ pairs in each block yields
		\[
		\begin{aligned}
			&\left\|\left(\sum_s\sum_{\substack{n\ge n_0,\ \ell\ge1\\s\le n+\ell<2s}}
			\sum_{j\ge a}\Gamma_j^s*g_{n-j,\ell}\right)_a
			\right\|_{L_2(\A;\ell_\infty)}
			\lesssim_d\|V\|_1(\alpha\|f\|_1)^{1/2},\\
			&\sum_s\sum_{\substack{n\ge n_0,\ \ell\ge1\\s\le n+\ell<2s}}
			\sum_j\|(H_j^s-\Gamma_j^s)*g_{n-j,\ell}\|_1
			\lesssim_d\|V\|_1\|f\|_1\sum_s s^22^{-s/8}
			\lesssim_d\|V\|_1\|f\|_1.
		\end{aligned}
		\]
		The first estimate uses $\sum_s s^22^{(2\delta_0-1/2)s}<\infty$.
		For the bounded good function, Lemma~\ref{lem:general-good-max} gives
		\[
		\|(S_a^Kg)_a\|_{L_2(\A;\ell_\infty)}^2
		\lesssim_d A_K^2\|g\|_2^2
		\lesssim_d A_K^2\alpha\|f\|_1.
		\]
		Apply Lemmas~\ref{lem:L2-cut} and~\ref{lem:absolute-L1} and intersect
		the resulting projections as in Section~\ref{sec:main-proof}.
		With $\alpha=\lambda/A_K$ and $\|V\|_1\le A_K$, this gives
		\begin{equation}\label{eq:general-smooth-endpoint}
			\sup_a\|eS_a^Kfe\|_\infty\le\lambda/2,\qquad
			\varphi(1-e)\lesssim_d
			\left[\alpha^{-1}+A_K\lambda^{-1}
			+A_K^2\alpha\lambda^{-2}\right]\|f\|_1
			\lesssim_d\frac{A_K\|f\|_1}{\lambda}.
		\end{equation}
		
		For $2^a\le\varepsilon<2^{a+1}$, we have 
		\[
		\left|K(y)\big(\one_{\{|y|>\varepsilon\}}
		-\rho(2^{-a}|y|)\big)\right|
		\le\sum_{j=a-1}^{a+1}|K_j(y)|.
		\]
		The same real and imaginary part compression as in
		Section~\ref{sec:main-proof} yields
		\[
		\|e(T_{K,\varepsilon}-S_a^K)fe\|_\infty
		\le2\sum_{j=a-1}^{a+1}\|e(|K_j|*f)e\|_\infty.
		\]
		Intersect the projection in \eqref{eq:general-smooth-endpoint} with
		that of Lemma~\ref{lem:positive-annuli} at height $\lambda/12$.
		This proves the asserted estimate  under \eqref{eq:finitemodel}.
		The general case follows as in Section~\ref{sec:main-proof}.
	 
	\end{proof}
	\section*{Acknowledgments}

	During the preparation of this work, the authors used ChatGPT   to assist with language editing and
	refinement of mathematical arguments. The authors take full responsibility for the content of this paper.

\end{document}